\documentclass[11pt,draftcls,onecolumn]{IEEEtran}

\usepackage[sort]{cite}
\usepackage{psfrag}
\usepackage{epsfig}
\usepackage{amsmath,bbm,times,stmaryrd}

\usepackage{amsthm}

 \usepackage[mathscr]{eucal}

\usepackage{tabularx}
\usepackage{multirow}
\usepackage{enumerate}

\usepackage[normal]{subfigure}

\usepackage{colortbl}
\usepackage{dcolumn}
\newcolumntype{.}{D{.}{.}{1.3}}

\usepackage{txfonts}
\usepackage[subnum]{cases}

\usepackage{equationarray}

\newcommand\dashrule{\leavevmode\xleaders\hbox{-}\hfill\kern0pt}

\def\argmin{\operatornamewithlimits{arg\,min}}

\newcommand{\bg}{\bvec{g}}

\newcommand{\bt}{\bvec{t}}
\newcommand{\bu}{\bvec{u}}
\newcommand{\bv}{\bvec{v}}
\newcommand{\bw}{\bvec{w}}
\newcommand{\bx}{\bvec{x}}
\newcommand{\by}{\bvec{y}}
\newcommand{\bz}{\bvec{z}}

\newcommand{\bA}{{\bf A}}
\newcommand{\bB}{{\bf B}}
\newcommand{\bC}{{\bf C}}
\newcommand{\bD}{{\bf D}}

\newcommand{\bF}{{\bf F}}

\newcommand{\bH}{{\bf H}}
\newcommand{\bI}{{\bf I}}

\newcommand{\bK}{{\bf K}}

\newcommand{\bP}{{\bf P}}
\newcommand{\bQ}{{\bf Q}}
\newcommand{\bR}{{\bf R}}

\newcommand{\bT}{{\bf T}}
\newcommand{\bU}{{\bf U}}
\newcommand{\bV}{{\bf V}}
\newcommand{\bW}{{\bf W}}
\newcommand{\bX}{{\bf X}}
\newcommand{\bY}{{\bf Y}}
\newcommand{\bZ}{{\bf Z}}

\newcommand{\calD}{{\mathcal{D}}}

\newcommand{\calL}{{\mathcal{L}}}

\newcommand{\bdelta}{\mbox{\boldmath $\delta$}}

\newcommand{\boldeta}{\mbox{\boldmath $\eta$}}
\newcommand{\btheta}{\mbox{\boldmath $\theta$}}

\newcommand{\bsigma}{\mbox{\boldmath $\sigma$}}

\newcommand{\1}{\mbox{\boldmath $1$}}
\newcommand{\0}{\mbox{\boldmath $0$}}

\newcommand{\be}{\begin{eqnarray}}
\newcommand{\ee}{\end{eqnarray}}

\newcommand{\matrixb}{\left[ \begin{array}}
\newcommand{\matrixe}{\end{array} \right]}

\def\*{\circledast}

\newcommand{\bvec}[1]{\boldsymbol{#1}}

\newcommand{\ve}{\bvec{e}}

\def\vectorize{\operatorname{vec}}
\newcommand{\vtr}[1]{\vectorize\hspace{-.3ex}\left(#1\right)}

\newcommand{\tensor}[1]{\boldsymbol{\mathscr{\MakeUppercase{#1}}}} 

\newcommand{\tE}{\tensor{E}}

\newcommand{\tI}{\tensor{I}}

\newcommand{\tW}{\tensor{W}}
\newcommand{\tX}{\tensor{X}}
\newcommand{\tY}{\tensor{Y}}

\usepackage[vlined,ruled,commentsnumbered]{algorithm2e}

\usepackage{aurical}

\def\blkdiag{\operatorname{blkdiag}}

\usepackage{arydshln}

\usepackage{amssymb}
\usepackage{amsmath}
\usepackage{amsfonts}
\usepackage{mathdots}

\usepackage{url}

\usepackage{mathtools}

\usepackage{graphicx}
\graphicspath{{./Linregbound/}{./Bestrank1/}{./PARO/}}

\definecolor{razzledazzlerose}{rgb}{1.0, 0.2, 0.8}
\definecolor{amethyst}{rgb}{0.6, 0.4, 0.8}
\definecolor{chromeyellow}{rgb}{1.0, 0.65, 0.0}

\title{Best Rank-One Tensor Approximation and Parallel Update Algorithm for CPD}
\author{Anh-Huy Phan, Petr Tichavsk\'y and Andrzej Cichocki
\thanks{A.-H. Phan and A. Cichocki are with the Lab for Advanced Brain Signal Processing, Brain Science Institute, RIKEN, Wakoshi, Japan, e-mail: (phan,cia)@brain.riken.jp.}
\thanks{A. Cichocki is also with Systems Research Institute PAS, Warsaw, Poland, and Skolkovo Institute of Science and Technology (Skoltech), Russia}
\thanks{P.  Tichavsk{\'y} is with Institute of Information Theory and Automation, Prague, Czech Republic, email: tichavsk@utia.cas.cz.}
\thanks{The work of P. Tichavsk{\'y} was supported by the Czech Science Foundation through project No. 17--00902S.}}

\renewcommand{\CommentSty}[1]{\fontsize{8.7}{9.8}\selectfont\textnormal{\texttt{#1}}\unskip}
\SetKwComment{mtcc}{\% }{}
\SetKwSwitch{Switch}{Case}{Other}{switch}{}{case}{otherwise}{end switch}%

\newcounter{example} 
    
\newenvironment{example}
{\refstepcounter{example}\vspace{10pt}\par\noindent 
\textbf{Example \theexample\ }
}
{}%

\makeatletter
\newsavebox{\@brx}
\newcommand{\llangle}[1][]{\savebox{\@brx}{\(\m@th{#1\langle}\)}%
  \mathopen{\copy\@brx\kern-0.5\wd\@brx\usebox{\@brx}}}
\newcommand{\rrangle}[1][]{\savebox{\@brx}{\(\m@th{#1\rangle}\)}%
  \mathclose{\copy\@brx\kern-0.5\wd\@brx\usebox{\@brx}}}
\makeatother

\emergencystretch=\maxdimen
\def\comment#1{}

\begin{document}
\sloppy

\maketitle

\begin{abstract}
A novel algorithm is proposed for CANDECOMP/PARAFAC tensor decomposition
to exploit best rank-1 tensor approximation. Different from the existing algorithms, our algorithm updates rank-1 tensors simultaneously in-parallel. In order to achieve this, we develop new all-at-once algorithms for best rank-1 tensor approximation based on the Levenberg-Marquardt method and the rotational update.  
We show that the LM algorithm has the same complexity of first-order optimisation algorithms, while the rotational method leads to solve the best rank-1 approximation of tensors of size $2 \times 2 \times \cdots \times 2$. We derive closed-form expression of best rank-1 tensor of $2\times 2 \times 2$ tensors, and present an ALS algorithm which updates 3 component at a time for higher order tensors.
The proposed algorithm is illustrated in decomposition of difficult tensors which are associated with multiplications of two matrices.
\end{abstract}

\section{Introduction}

The CANDECOMP/PARAFAC tensor decomposition seeks the best rank-$R$ tensor approximation to a data tensor $\tY$ of size $I_1 \times I_2 \times \cdots \times I_N$
\be
	\tY \approx  \sum_{r = 1}^{R} \tX_r \notag 
\ee
where $	\tX_r = \bu_{1,r} \circ \bu_{2,r} \circ \cdots \circ \bu_{N,r}$ are rank-1 tensors.
Matrix of the loading components,  $\bU_{n} = [\bu_{n,1}, \ldots, \bu_{n,R}]$, are coined the factor matrices.


%
%


The CPD can be achieved by minimising the Frobenius norm of the error 
\be
\min \quad \|\tY - \sum_{r = 1}^{R} \tX_r \|_F^2	\label{eq_cpd_cost} \, .
\ee
The above objective function is nonlinear with respect to all the factor matrices, but linear in parameters in one factor matrix, 
or parameters in non-overlapping partitions of different factor matrices \cite{Petr_PALS,Phan_PHALS}. Hence, one can derive the Alternating Least Squares (ALS) to update sequentially the factor matrices $\bU_n$\cite{Harshman,Phan_fastALS}, or jointly update the loading components in non overlapping partitions \cite{Petr_PALS,Phan_PHALS}. An alternative method is to apply the  nonlinear conjugate gradient method \cite{Paatero94,JChem-CPOPT}, the Levenberg-Marquardt algorithm \cite{Paatero97,Phan_fLM}, the non-linear least squares  (NLS) algorithm\cite{Sorber_fusion2013} to update all the parameters at a time. 

Consider a particular case of CPD with rank-1, i.e., finding a best rank-1 tensor approximation. This is the case where the Tucker decomposition (TKD), the tensor network, tensor train, hierarchical Tucker  decomposition meet CPD. Following this, the best rank-1 tensor approximation inherits good algorithms from the other tensor network decompositions, such as the sequential projection and truncation method, also known TT-SVD, or the DMRG algorithm for Tensor-train tensor decomposition. 
The sequential projection and truncation has been recently shown to be a good method for finding the best rank-1 tensor approximation \cite{TruncatedSVD2011,da2016finite}.

In addition, the well-known Higher Order Orthogonal Iteration (HOOI) algorithm for TKD becomes the ALS algorithm for CPD. 

The best rank-1 tensor approximation is the only one case when the Levenberg-Marquardt (LM) algorithm has a similar computational cost to that of ALS. Its update rule can be proved to be in a similar form of the first order optimization algorithm with an optimally determined step-size. 

In addition, all the loading components can be updated through a rotational method, which in turn solves a best rank-1 tensor approximation to a quantised-scale tensor of size $2 \times 2 \times \cdots \times 2$.

With the good algorithms for best rank-1 tensor approximation, the questions are   

\begin{itemize}
\item ``\emph{Can we employ the best rank-1 tensor approximation for higher rank CPD?}''

\item ``\emph{Can we update rank-1 tensors in CPD simultaneously in parallel?}''
\end{itemize}

Indeed the high rank CPD can be formulated as a sequence of best rank-1 tensor approximations to the residue tensors 
\be
	\min \quad \|\tY_r - \bu_{1,r} \circ \bu_{2,r} \circ \cdots \circ \bu_{N,r} \|_F^2 \label{eq_hals_obj}
\ee
where $\tY_r = \tY - \sum_{s \neq r} \tX_s  = \tE + \tX_r$.
This is the way to derive the hierarchical ALS (HALS) algorithm \cite{Phan-HALS-IEICE}.
 Once a rank-1 tensor $\tX_r$ is updated, HALS updates the error tensor $\tE$ and proceeds the next rank-1 tensor approximation. Since the error tensor $\tE$ varies in  sub-problems, HALS cannot update all rank-1 tensors, $\tX_1$, \ldots, $\tX_R$, simultaneously.

So far, tensor deflation is the only one method able to extract rank-1 tensors in parallel \cite{Phan_tensordeflation_alg,Phan_ranksplitting_full}. However, this kind of tensor decomposition requires additional conditions, and does not rely on the best rank-1 tensor approximation.

In this paper, we address the above two questions, and propose a novel algorithm which can update rank-1 tensors simultaneously in parallel. Moreover, we derive a novel LM algorithm and a rotational algorithm for the best rank-1 tensor approximation. These proposed algorithms have the same computational cost as that of the ALS/HOOI algorithm.

\section{Parallel Rank-1 Tensor Update Algorithm}

We denote vectorize of the data tensor $\tY$ by $\by = \vtr{\tY}$, and a matrix $\bX$ comprising vectorisation of rank-1 tensors $\tX_r$
\be
\bX &=& [\vtr{\tX_1}, \vtr{\tX_2}, \ldots, \vtr{\tX_R}]\, \notag \\
   &=&  \bU_N \odot \cdots \odot \bU_2 \odot \bU_1 \label{eq_X_mat}\,.
\ee
$\bX$ is also known as Khatri-Rao product of the factor matrices $\bU_1, \ldots, \bU_N$.
Now, we rephrase the CPD in (\ref{eq_cpd_cost}) in a new form to find a matrix $\bX$ holding the Khatri-Rao structure, that is, 
\be
\min \quad && f(\bX) = \frac{1}{2} \, \|\by - \bX\,  \1_R\|_2^2   \label{eq_cpd_cost_2} \\
\text{s.t.} \quad&&   \bX = \bU_N \odot \cdots \odot \bU_2 \odot \bU_1\notag\,,
\ee
where $\1_R$ is a vector of ones of the length $R$. This can be interpreted as a generalised projection problem 
\be
\min  \quad && f(\bZ) + g(\bX)  \label{eq_cpd_cost_3} \, , \quad \text{s.t.} \quad \bZ = \bX \, , \notag 
\ee
where $g(\bX)$ is the indicator function of a set, $\calD$, of structured matrices which are in form of the Khatri-Rao products (\ref{eq_X_mat}), i.e., $g(\bX) = 0$ if $\bX \in \calD$, otherwise $\infty$. The problem can be solved using the augmented Lagrangian, or more specifically the alternating direction method of multipliers (ADMM) or the alternating projection method  \cite{Proximal2014,Becker2011}. The augmented Lagrangian function to the problem (\ref{eq_cpd_cost_3}) is given by 
\be
	\calL(\bX, \bZ, \bT)  = f(\bZ) + g(\bX)  -  \frac{1}{\gamma} \,   \langle \bZ - \bX, \bT  \rangle + \frac{1}{2\gamma} \, \|\bZ - \bX\|_F^2    \label{eq_lagrangian}\, ,
\ee
where $\gamma > 0$, $\bZ$ and $\bX$ are primal variables, and $\bT$ is the dual variable, $\langle \bA, \bB \rangle$ denotes the scalar product of two matrices.
The Lagrangian function can be rewritten in the form of 
\be
	\calL(\bX, \bZ, \bT)  = f(\bZ) + g(\bX) + \frac{1}{2\gamma} \, \left(\|\bZ - \bX - \bT \|_F^2 - \|\bT\|_F^2\right)  \,.\label{eq_lagrangian2}
\ee
Updates of the primal variables, $\bZ$, $\bX$, and the dual variable, $\bT$, consist of the following iterations  
\begin{align}
\bZ^{(k+1)} &= \argmin_{\bZ} \;  \frac{1}{2} \, \|\by - \bZ \, \1_R \|_F^2 + \frac{1}{2\gamma} \, \|\bZ - \bX^{(k)} - \bT^{(k)} \|_F^2  \,,\label{eq_updateZ}\\
\bX^{(k+1)} &=  \argmin_{\bX} \; g(\bX) +\frac{1}{2\gamma} \, \|\bZ^{(k+1)}  - \bT^{(k)} - \bX \|_F^2 \, \notag  \\ 
&=   \Pi_{\calD}(\bZ^{(k+1)} - \bT^{(k)})  \, ,\label{eq_updateX}\\
\bT^{(k+1)} &=   \bT^{(k)}  + \bX^{(k+1)} - \bZ^{(k+1)} \label{eq_update_T}\,,
\end{align}
where $k$ denotes the iteration index, and $\Pi_{\calD}(\bZ)$  the projection of $\bZ$ onto $\calD$.

%
%

\subsection{Update of $\bZ$}
Since the problem (\ref{eq_updateZ}) is quadratic, 
$\bZ$ is found in closed-form as 
\begin{align}
\bZ^{(k+1)} &= \argmin_{\bZ} \quad  \frac{1}{2} \, \|\by - \bZ \, \1_R \|_F^2 + \frac{1}{2\gamma} \, \|\bZ - \bX^{(k)} - \bT^{(k)} \|_F^2 \notag \\
&=  \left(\bX^{(k)} + \bT^{(k)} +  \gamma \, \by \, \1_R^T \right) 
\left(\bI +   \gamma \1_R \1_R^T \right)^{-1}   \notag \\
&=  \left(\bX^{(k)} + \bT^{(k)} +  \gamma \, \by \, \1_R^T \right) 
\left(\bI  -   \frac{\gamma}{1 + \gamma R} \1_R \1_R^T \right)\notag \\
&=  \frac{\gamma}{1 + \gamma R} \,  \by \,  \1_R^T + \left(\bX^{(k)} + \bT^{(k)}\right) \, \left(\bI  -   \frac{\gamma}{1 + \gamma R} \1_R \1_R^T \right) \label{eq_updateZ2} \, .
\end{align}

\subsection{Update of $\bX$}
From (\ref{eq_updateX}), 
columns of the Khatri-Rao matrix $\bX^{(k+1)}$ can be updated independently as best rank-1 tensor approximation to the tensors whose vectorizations are $(\bz_r^{(k+1)} - \bt_r^{(k)})$, i.e., 
\be
\min_{\displaystyle \{\bu_{n,r}\}} \quad  \frac{1}{2} \|\bz_r^{(k+1)} - \bt_r^{(k)} -  \bu_{N,r} \otimes \cdots \otimes \bu_{2,r} \otimes \bu_{1,r} \|^2  \,.\label{eq_update_un_r} 
\ee
We can apply the ALS/HOOI algorithm or the fLM algorithm \cite{Phan_fLM} to solve the above problem. More computationally efficient algorithms are presented in Section~\ref{sec::bestrank1}.

With the three update rules (\ref{eq_update_T}), (\ref{eq_updateZ2}) and (\ref{eq_update_un_r}), we can implement an algorithm to sequentially update $\bZ$, $\bX$ and $\bT$. However, such a simple algorithm demands a large extra space for the matrices $\bZ$ and $\bT$ of size $(I_1 \ldots I_N) \times R$, and even the matrix of rank-1 tensors, $\bX$. In the following subsection we present a memory saving implementation of the proposed procedure, which only requires
memory of the order $5 I_1 I_2 \cdots I_N$.

\subsection{Rank-1 Tensor Update in Parallel}

First, we consider the term $\bZ^{(k+1)} - \bT^{(k)}$ which appears in the update  (\ref{eq_update_un_r}) for $\bU_{n}$. 
%
From the update for $\bZ$ in (\ref{eq_updateZ2}), we have
\be
\bZ^{(k+1)} - \bT^{(k)}  =  \bX^{(k)} + \frac{\gamma R}{1 + \gamma R} (\bar{\by} -  \bar{\bx}^{(k)} - \bar{\bt}^{(k)}) \, \1_R^T\,  \label{eq_updateZ3}
\ee
where $\bar{\by} = \displaystyle \frac{1}{R} \by$, and $\displaystyle \bar{\bx}^{(k)} = \frac{1}{R} \, \bX^{(k)} \, \1_R$ and $\displaystyle \bar{\bt}^{(k)} = \frac{1}{R} \,  \bT^{(k)} \, \1_R$  are means of columns of $\bX^{(k)}$ and $\bT^{(k)}$, respectively. 
Note that $\bar{\bx}^{(k)}$ is the rank-$R$ tensor approximation of the tensor $\frac{1}{R} \tY$ but in the vectorisation form of $\bar{\by}$ at the iteration-$k$.

We define a parameter $\mu$ which depends on $\gamma$ and $R$ as
\be
\mu =  \frac{\gamma R}{1 + \gamma R} \, \label{eq_mu_paro}
\ee
and a residue at the iteration-$k$
\be
\ve^{(k)}  = \mu\, (\bar{\by} -  \bar{\bx}^{(k)} - \bar{\bt}^{(k)}) \,  \, \label{eq_residual_ek}.
\ee 
It is obvious from (\ref{eq_updateZ3}) that 
\be
	\bz_r^{(k+1)} - \bt_r^{(k)} = \bx_r^{(k)} +  \ve^{(k)} \, .\notag 
\ee
When the algorithm converges, the components $\bz_r$ are (nearly) in the form of the Kronecker product of the loading components $\bu_{n,r}$, and the dual variables $\bt_r$ become zeros or take small values. It follows that the residue $\ve$ will also go to zero.

The best rank-1 tensor approximation of $(\bz_r^{(k+1)} - \bt_r^{(k)})$ in (\ref{eq_update_un_r}) is rewritten as  
\be
\min_{\bu_{n,r}} \quad  \frac{1}{2} \|\ve^{(k)} +  \bx_r^{(k)} -  \bu_{N,r} \otimes \cdots \otimes \bu_{2,r} \otimes \bu_{1,r} \|^2  \label{eq_update_un_r2} \,.
\ee
Next, we replace $\bZ^{(k+1)}$ in (\ref{eq_updateZ2}) into the update of $\bT$ in (\ref{eq_update_T}) 
\be
\bT^{(k+1)} &=&  \bT^{(k)}  + \bX^{(k+1)}  - \left(\bX^{(k)} + \bT^{(k)}\right) \, \left(\bI  -   \frac{\gamma}{1 + \gamma R} \1_R \1_R^T \right)  \notag \\ && \quad -   \frac{\gamma}{1 + \gamma R} \,  \by \,  \1_R^T \notag \\
&=&    \bX^{(k+1)}  -  \bX^{(k)}   - \frac{\gamma R}{1 + \gamma R} \,    \left(\bar{\by} - \bar{\bx}^{(k)} -  \bar{\bt}^{(k)} \right) \,  \1_R^T \notag \\
&=&  \bX^{(k+1)}  -  \bX^{(k)}  -  \ve^{(k)}  \, \1_R^T \,.
\label{eq_update_t2}
\ee
From (\ref{eq_residual_ek}) and (\ref{eq_update_un_r2}), it reveals that there no need to compute the dual variables $\bT^{(k)}$, but only the term $\bar{\bt}^{(k)}$, which, from (\ref{eq_update_t2}) and its definition, is given by 
\be
\bar{\bt}^{(k+1)}  = \frac{1}{R} \bT^{(k+1)} \, \1_R &=&     \bar{\bx}^{(k+1)}  -  \bar{\bx}^{(k)}   -  \ve^{(k)} \, .\label{eq_update_t3}
\ee
Now we can even omit $\bar{\bt}^{(k)}$ by replacing its expression into (\ref{eq_residual_ek}) to obtain 
\be
\ve^{(k)}  = \mu\, (\bar{\by} -  2 \, \bar{\bx}^{(k)} +  \bar{\bx}^{(k-1)}  +  \ve^{(k-1)}) \,  \,. \label{eq_residual_ek2} 
\ee 
The final update of the residue $\ve$ does not comprise $\bar{\bt}$, but relates to the current error, $\bar{\by} - \bar{\bx}^{(k)}$,  and the error between two estimated tensors $(\ve^{(k-1)} + \bar{\bx}^{(k-1)})- \bar{\bx}^{(k)}$. More precisely, the residue $\ve$ is updated from the 2nd-order difference of the sequence $\bar{\bx}^{(k-1)} + \ve^{(k-1)}, \bar{\bx}^{(k)}, \bar{\by}$.

%

 \setlength{\algomargin}{1em}
\begin{algorithm}[t!]
\SetFillComment
\SetSideCommentRight
\CommentSty{\footnotesize}
\caption{{\tt{PArallel Rank-One Update (PARO)}}\label{alg_paro}}
\DontPrintSemicolon \SetFillComment \SetSideCommentRight
\KwIn{Data tensor $\tY$:  $(I_1 \times I_2 \times \cdots \times I_N)$,  and a rank $R$ and a regularisation parameter $\mu$ in $[\frac{1}{2}, 1]$} 
\KwOut{$\tX = \llbracket  \bU^{(1)}, \bU^{(2)}, \ldots, \bU^{(N)}\rrbracket $ of rank $R$}
\SetKwFunction{cpd}{CPD} 
\SetKwFunction{bals}{bals} 
\SetKwFor{For}{parfor}{}{}
\Begin{
\nl Initialize $\tX^{(0)} = \llbracket  \bU_1^{(0)}, \bU_2^{(0)}, \ldots, \bU_N^{(0)}\rrbracket$ \;
\nl $\bar{\by} = \frac{1}{R} \vtr{\tY}$\;
\Repeat{a stopping criterion is met}{
\nl $\bar{\bx}^{(k)} = \frac{1}{R} \sum_{r}  {\bx}_r^{(k)}$\;
\nl $\ve^{(k)}  = \mu \, (\bar{\by} - 2 \bar{\bx}^{(k)} + \bar{\bx}^{(k-1)}  +  \bar{\ve}^{(k-1)}) $\;
 \For(\tcc*[f]{process in parallel}){$r = 1, 2, \ldots, R$}
{
\nl Find ${\bx}_r^{(k+1)}$ as the best rank-1 tensor approximation of $\ve^{(k)} +  \bx_r^{(k)} $ \\  $\min_{\bu_{n,r}} \quad  \frac{1}{2} \|\ve^{(k)} +  \bx_r^{(k)} -  \bu_{N,r} \otimes \cdots \otimes \bu_{2,r} \otimes \bu_{1,r} \|^2 $} \label{step_update_rank1}
}
}
\end{algorithm}

\subsection{Implementation}\label{sec:implement_cp_bound}
 
Algorithm~\ref{alg_paro} shows a simple  implementation of the proposed PArallel Rank-One tensor update algorithm (PARO) which consists of the two update rules in (\ref{eq_update_un_r2}) and  (\ref{eq_residual_ek2})
\begin{itemize}
\item Update the residue $\ve^{(k)}$ as in (\ref{eq_residual_ek2})
\item Seek in parallel best rank-1 tensors $\tX_r^{(k+1)}$ to the residue tensors whose vectorisations are $\ve^{(k)} + \bx_r^{(k)}$ for $r = 1, 2, \ldots, R$. This step is performed in a distributed system with multi-nodes, each node estimates one rank-1 tensor.

\end{itemize}

\subsubsection{Similarity with HALS}

In the above implementation, PARO works in a similar way to the HALS algorithm which minimises the objective function in (\ref{eq_hals_obj}) \cite{Phan-HALS-IEICE}. However, the residue in HALS is the error between the data tensor $\tY$ and the current estimate $\tX^{(k)}$, i.e., $\tE^{(k)} = \tY - \tX^{(k)}$ and $\tY_r = \tE^{(k)} + \tX_r^{(k)}$, and it is updated after each rank-1 tensor approximation, while PARO uses the same residue tensor in all $R$ rank-1 tensor approximations. 

\subsubsection{Construction of rank-1 tensor components $\tX_r$}

In our algorithm, the rank-1 tensors, $\tX_r$, or their vectorizations, $\bx_r$, appear in the sub-problems in (\ref{eq_update_un_r2}) and in the update of $\bar{\bx}^{(k)}$. However, we need not save all $R$ rank-1 tensors. Section~\ref{sec::bestrank1} will show that 
in the ALS, LM and rotational update algorithms for best rank-1 tensor approximation, the tensors $\tE^{(k)} + \tX_r^{(k)}$ involve only in the following tensor-vector products, e.g., in (\ref{eq_update_r1als}), (\ref{eq_tn}) 
\be	
&&(\tE^{(k)} + \tX_r^{(k)})  \, \bar{\times}_1 \, \bu_{r,1}  \, \bar{\times}_2 \, \bu_{r,2}  \cdots \bar{\times}_N \, \bu_{r,N}  \notag \\&&\quad= \tE^{(k)} \, \bar{\times}_1 \, \bu_{r,1}  \, \bar{\times}_2 \, \bu_{r,2}  \cdots \bar{\times}_R \, \bu_{r,R} +  \prod_{n = 1}^{N} \, (\bu_{r,n}^T \, \bu_{r,n}^{(k)})  \notag 
\ee
or the tensor-matrix products in (\ref{eq_f_eta}) and (\ref{eq_W_proj})
\be	
&&(\tE^{(k)} + \tX_r^{(k)})  \, \bar{\times}_1 \, \bV_1\,   \cdots \bar{\times}_N \,  \bV_N 
\notag \\
&& \quad = \tE^{(k)} \, \bar{\times}_1 \, \bV_1 \,   \cdots \bar{\times}_N \,  \bV_N+  (\bV_1^T \bu_{1,r} ) \circ \cdots  \circ (\bV_N^T \bu_{N,r} )  \notag 
\ee 
where $\bV_n= [\bu_{r,n},  \bg_{n}]$ consist of two columns, and $\bY \bar{\times}_n \bV = \bY  {\times}_n \bV^T$. 
The first case is related to scalar products, while for the latter  case, the products $(\bV_N^T \, \bu_{N,r} )$ result vectors of length 2.
Hence, there is no need to construct explicitly the tensors $\tE^{(k)} + \tX_r^{(k)}$.

Following steps listed in Algorithm~\ref{alg_paro}, the proposed algorithm needs $2 I_1 I_2 \cdots I_N$ memory cells to compute $\bar{\bx}^{(k)}$, space to store $\bar{\bx}^{(k-1)}$ and $\ve^{(k)}$. In total, it needs a space for $5 I_1 I_2 \cdots I_N$ entries for $\bar{\bx}^{(k-1)}$, $\bar{\bx}^{(k)}$, $\ve^{(k)}$, $\bar{\by}$ and a temporary parameter.

\subsubsection{Choice of the regularisation parameter $\mu$}

Another important factor in our augmented Lagrangian based algorithm is the choice of the regularisation parameter $\gamma$ or the parameter $\mu$ defined in (\ref{eq_mu_paro}).
Similar to the alternating direction method of multipliers and the generalised projection method, the algorithm may diverge with an unsuitable step size \cite{BoydFTML2011,Becker2011}.

At the beginning of the estimation process, $\ve^{(0)} = \bar{\by} - \bar{\bx}^{(0)}$, hence $\ve^{(1)} = 2\mu (\bar{\by} - \bar{\bx}^{(1)})$, and 
we can choose
\be
\mu = \frac{1}{2}, \quad \text{i.e.,} \quad  \gamma = \frac{1}{R} \,. \notag 
\ee
This natural choice of $\mu$ keeps $\ve = \bar{\by}  - \bar{\bx}$ and is considered a default value in PARO.
However, it may not be the best, and can make the PARO algorithm converge slowly. 
We present efficient strategies to select and adjust $\gamma$ or $\mu$, and illustrate them through Example~\ref{ex_nasob222}.

\begin{example}[Effect of regularisation parameters on the convergence of PARO.]\label{ex_nasob222}
We illustrate performance of PARO for decomposition of the tensor for multiplication of two matrices of size $2 \times 2$. 
This tensor is of size $4 \times 4 \times 4$, contains only zeros and ones, and its four frontal slices are given by
\be
\bY_1 =   \left[\arraycolsep=1.6pt\def\arraystretch{.6}
\begin{array}{cccc}	 
1 &    0 &    0&     0\\
0 &    0  &   1 &    0\\
0 &    0  &   0 &    0\\
0 &    0  &   0 &    0
\end{array}
\right],   \, 
\bY_2 =
\left[\arraycolsep=1.6pt\def\arraystretch{.6}
\begin{array}{cccc}	 
0 &    0 &    0&     0\\
0 &    0  &   0 &    0\\
1 &    0  &   0 &    0\\
0 &    0  &   1 &    0
\end{array}
\right],  \,   
\bY_3 = \left[
\arraycolsep=1.6pt\def\arraystretch{.6}
\begin{array}{cccc}	 
0 &    1 &    0&     0\\
0 &    0  &   0 &    1\\
0 &    0  &   0 &    0\\
0 &    0  &  0 &    0
\end{array}
\right], \,
\bY_4 =
\left[
\arraycolsep=1.6pt\def\arraystretch{.6}
\begin{array}{cccc}	 
0 &    0 &    0&     0\\
0 &    0  &   0 &    0\\
0 &    1  &   0 &    0\\
0 &    0  &  0 &    1
\end{array}
\right]\,, \notag 
\ee 
 and obey 
\be
	\vtr{\bA \bB} = \tY \times_1 \vtr{\bA^T}^T \times_2 \vtr{\bB^T}^T \, \notag
\ee
for any matrices $\bA$ and $\bB$ of the size $2 \times 2$.
The tensor is considered of rank-$R=7$ \cite{Strassen:1969:GEO:2722431.2722798,journals/corr/TichavskyPC16}.

In spite of a relatively small  tensor, decomposition of this tensor using the ordinary ALS requires a thousands of iterations as shown in Fig.~\ref{fig_222_1}.

{\bf{PARO with a fixed regularisation parameter.}}

We can run PARO with a fixed regularisation parameter, 
e.g., $\gamma R = 1$, or try several higher values of $\gamma R$, then choose the value which gives a good convergence.

Using the same initialization as ALS and with the error preservation norm correction of the initial \cite{Phan_EPC}, PARO with a default setting of $\mu = \frac{1}{2}$ converged after at most 2500 iterations. 

\begin{figure}[t!]
\centering
\subfigure[ALS vs PARO with various $\mu$.]{\includegraphics[width=.48\linewidth, trim = 0cm  0cm 0cm 0cm,clip=true]{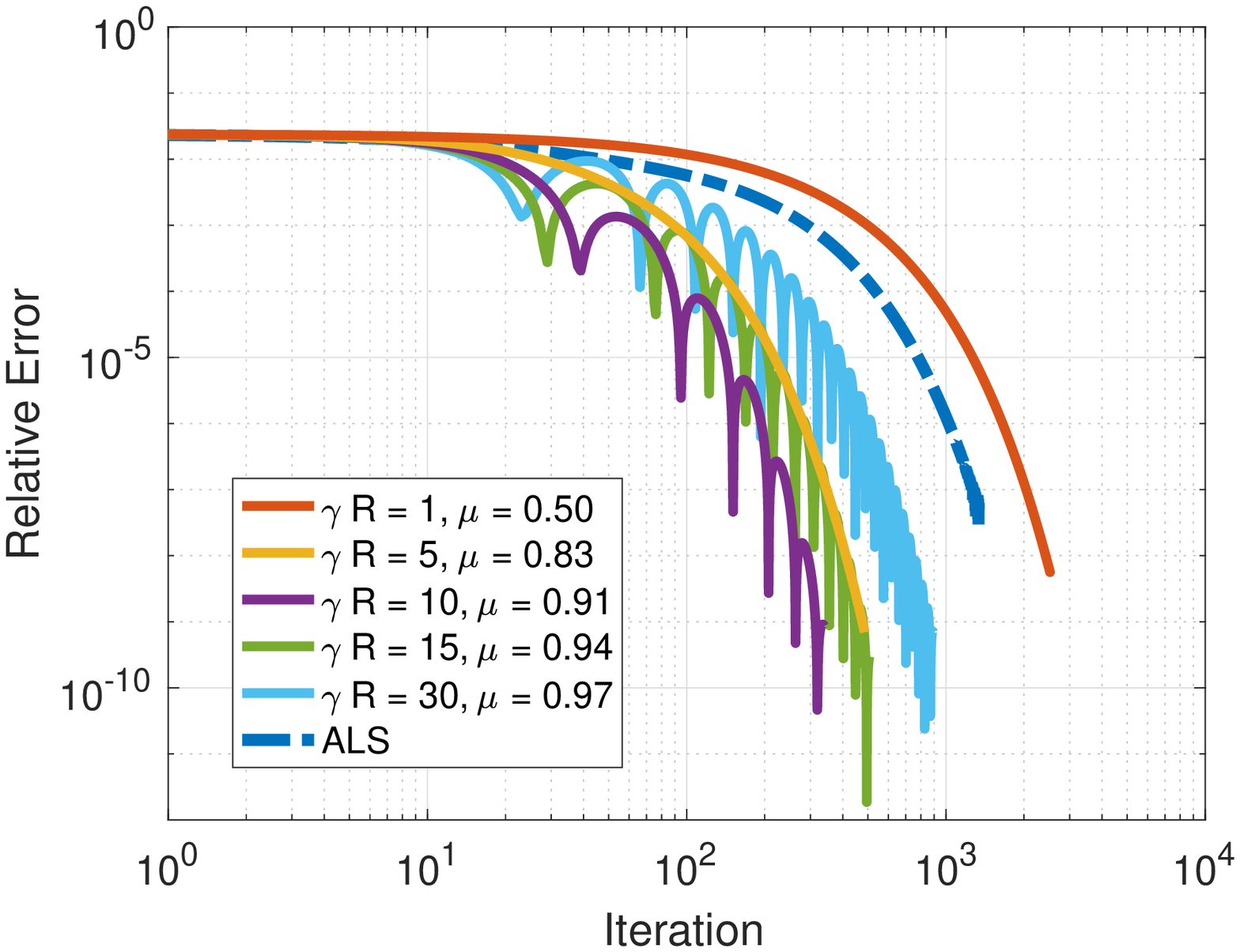}\label{fig_222_1}}
\subfigure[PARO with various adjusting rates.]{\includegraphics[width=.48\linewidth, trim = 0cm 0cm 0cm 0cm,clip=true]{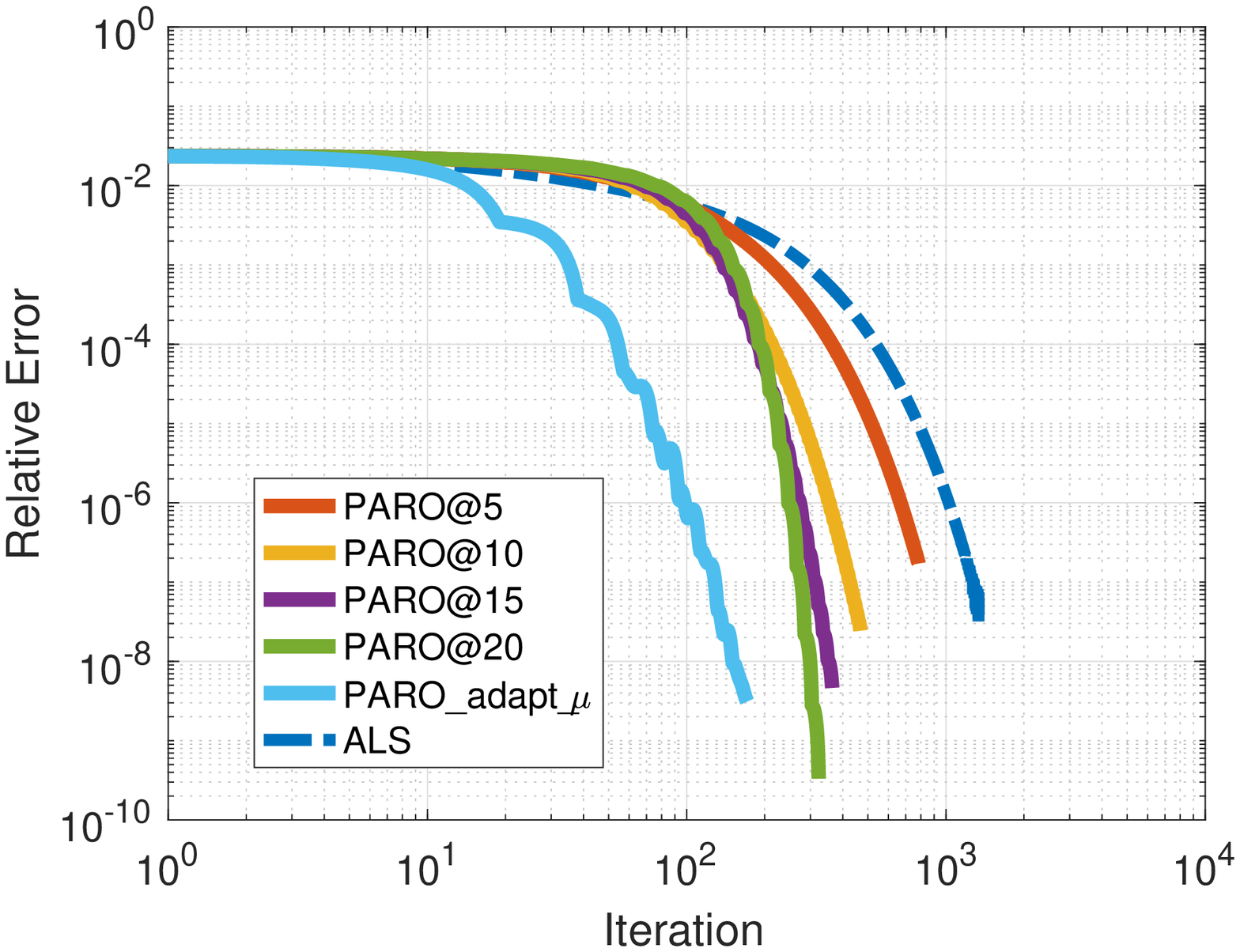}\label{fig_222_2}
}
\caption{Convergence of PARO with various settings of the regularisation parameter $\gamma$ or $\mu$ in Example~\ref{ex_nasob222}. \subref{fig_222_1} PARO with a fixed regularisation parameter $\gamma$. \subref{fig_222_2} $\gamma$ is adjusted by a factor of $\eta = \sqrt{2}$ every 5, 10, \ldots iterations or adaptively corresponding to the change of the objective values. }
\label{fig_222}
\end{figure}

A higher value of $\mu$, which is closer to 1, tends to decrease the objective function quickly, e.g., $\mu = 0.83$ corresponding to $\gamma R = 5$.
However, a relatively higher value of $\mu$,  e.g., $\mu = 0.94, 0.97$ as seen in Fig.~\ref{fig_222_1}, may make the algorithm unstable after a dozen of iterations. 
Similar behaviour was observed in other generalised projection method as discussed in \cite{Becker2011}. 
Despite of that, in this example, the algorithm still converges.

In this example, $\gamma R = 5$ is a good choice, and PARO converges in 486 iterations, faster than ALS with 1340 iterations. The algorithm can even converge faster with $\gamma R = 10$, but the objective function does not always decrease. 
\end{example}

{\bf{PARO with a regularly adjusted regularisation parameter.}}

Instead of specifying a fixed regularisation parameter, we can adaptively adjust it during the estimation process, e.g., 
decrease $\gamma$ every 10, 20 iterations, if the objective function is non-decreasing, otherwise, increase it.
A simple but efficient strategy is that we first execute PARO with $\gamma  = \frac{1}{R}$, i.e., $\mu = \frac{1}{2}$, then adjust $\gamma$ by a factor $\eta > 1$ every 10, 20 or 30 iterations, i.e.,
\be
\gamma \leftarrow  \,  \eta \, \gamma  , \quad \text{or}\quad 
\mu = \frac{1}{1 + \eta^{-k}}\,.\notag 
\ee
The parameter $\gamma$ might decrease by the same factor if the current setting $\gamma$ does not keep the objective function non-increasing.
Our experience is that $\eta = \sqrt{2}$ works in most experiments, including decompositions of the multiplication tensors and synthetic tensors which admit the considered model.
Discussion on step size adaptation can be further found in Sections 5.3 and 5.5 in \cite{Becker2011}.

For Example~\ref{ex_nasob222}, this update strategy speeds up the convergence of PARO as illustrated in Fig.~\ref{fig_222_2}. The results also indicate that increasing $\gamma$ or $\mu$ too fast might not improve much the convergence of PARO, although the algorithm is still able to converge faster than ALS. 
In addition, Fig.~\ref{fig_222_2} illustrates convergence of PARO with an adaptive adjustment of $\gamma$ for every 20 iterations, with an initial value of $5/R$.

\section{Best Rank-1 Tensor Approximation, where Tucker Decomposition Meets CPD}\label{sec::bestrank1}

We next present efficient algorithms for the best rank-1 tensor approximation, which are employed in the update in (\ref{eq_update_un_r2}) and Step~\ref{step_update_rank1} in Algorithm~\ref{alg_paro}. 
The section begins with a simple ALS algorithm and presents an efficient initialisation based on sequential truncation and projection. We show that the LM algorithm for this particular tensor decomposition has an equivalent form using the first order optimisation method. Finally, we propose a rotational algorithm for the best rank-1 tensor approximation.

\subsection{HOOI and HALS}
The simplest case of the tensor decomposition is with rank-1, i.e., seeking the best rank-1 tensor approximation
\be
\min \quad \|\tY - \bu_{1} \circ \bu_{2}  \circ \cdots \circ \bu_{N} \|_F^2 \, . \label{eq_bestrank1}
\ee
This is the case where the Tucker tensor decomposition (TKD) and tensor train meet  CPD.
The state-of-the-art Higher Order Orthogonal Iteration (HOOI) algorithm for TKD \cite{Lathauwer_HOOI} estimates loading components as principal components of symmetric matrices which in this case are of rank-1. The algorithm becomes the ALS or Hierarchical ALS algorithm \cite{Phan-HALS-IEICE}, where the loading components are projected vectors of the data $\tY$ by $(N-1)$ other loading components
\be
\bu_n =  \frac{\tY \, \bar{\times}_{k \neq n} \, \bu_k}{\displaystyle\prod_{k \neq n} (\bu_k^T \bu_k)}\, , \notag 
\ee
or 
\be
\bu_n =  \tY \, \bar{\times}_{k \neq n} \, \bu_k\,   \label{eq_update_r1als}
\ee
provided that the loading components are $\ell_2$ unit length vectors.
The estimated components are then normalised before proceeding the next iteration to update $\bu_{n+1}$.
It is know that such alternating algorithms can get stuck into local minima,
and its convergence depends on the initial values. A common method is to use leading left singular vectors of mode-$n$ matricization of the tensor, $\bY_{(n)}$, and 
more efficiently the sequential projection and truncation as in the TT-SVD algorithm for the Tensor train \cite{Vidal03,OseledetsTT09} (see next section and Example~\ref{ex_svd_tt_svd} for comparison between the two methods). 

\subsection{TT-SVD, the sequential projection and truncation method}

The TT-SVD \cite{Vidal03,OseledetsTT09,OseledetsTT11} 
was developed for the tensor train decomposition, in which core tensors of the tensor network are of order-2 or 3 and interconnected. When all TT-ranks are 1, TT-SVD serves for  best rank-1 tensor approximation. 
The first loading component $\bu_1$ is the leading left singular vector of the reshaping matrix $\bY_1 = \bY_{(1)}$
\be
	\bY_{1} \approx \bsigma_1 \, \bu_1 \, \bv_1^T   \notag \label{eq_tt_truncation} \, .
\ee
The right singular vector or the projected data $\sigma_1 \bv_1$ is then reshaped into a matrix $\bY_2$ of size $I_2 \times (I_3I_4\cdots I_N)$, and the second loading component $\bu_2$ is the leading left singular vector of this matrix.
The algorithm
executes $(N-1)$ sequential data projections and truncated-SVD in order to find $N$ loading components. The obtained components are then used to initialise the algorithms for best-rank-1 tensor approximation, e.g., HALS \cite{Phan-HALS-IEICE}, the Alternating Single or  Double-Core Update \cite{Phan_TT_part1}.

Since TT-SVD performs the data projections sequentially over  modes of the tensor, different combination of the tensor modes in the order of the tensor projection may lead to different results. Some of them are even worse than that using the SVD-based initialisation method.   
In other words, performance of TT-SVD highly depends on the projection order of the tensor modes. We can apply TT-SVD to various permutations of the tensor, then choose the best  result. 
In total, there are $N!$ combination of tensor modes. 
For tensors of low order, e.g., order-3, 4, we can run the algorithm for 6 or 24 tensor permutations.

\begin{figure}[t!]
\centering
\subfigure[$N = 3$]{\includegraphics[width=.48\linewidth, trim = 1.0cm 1cm 0cm 0cm,clip=true]{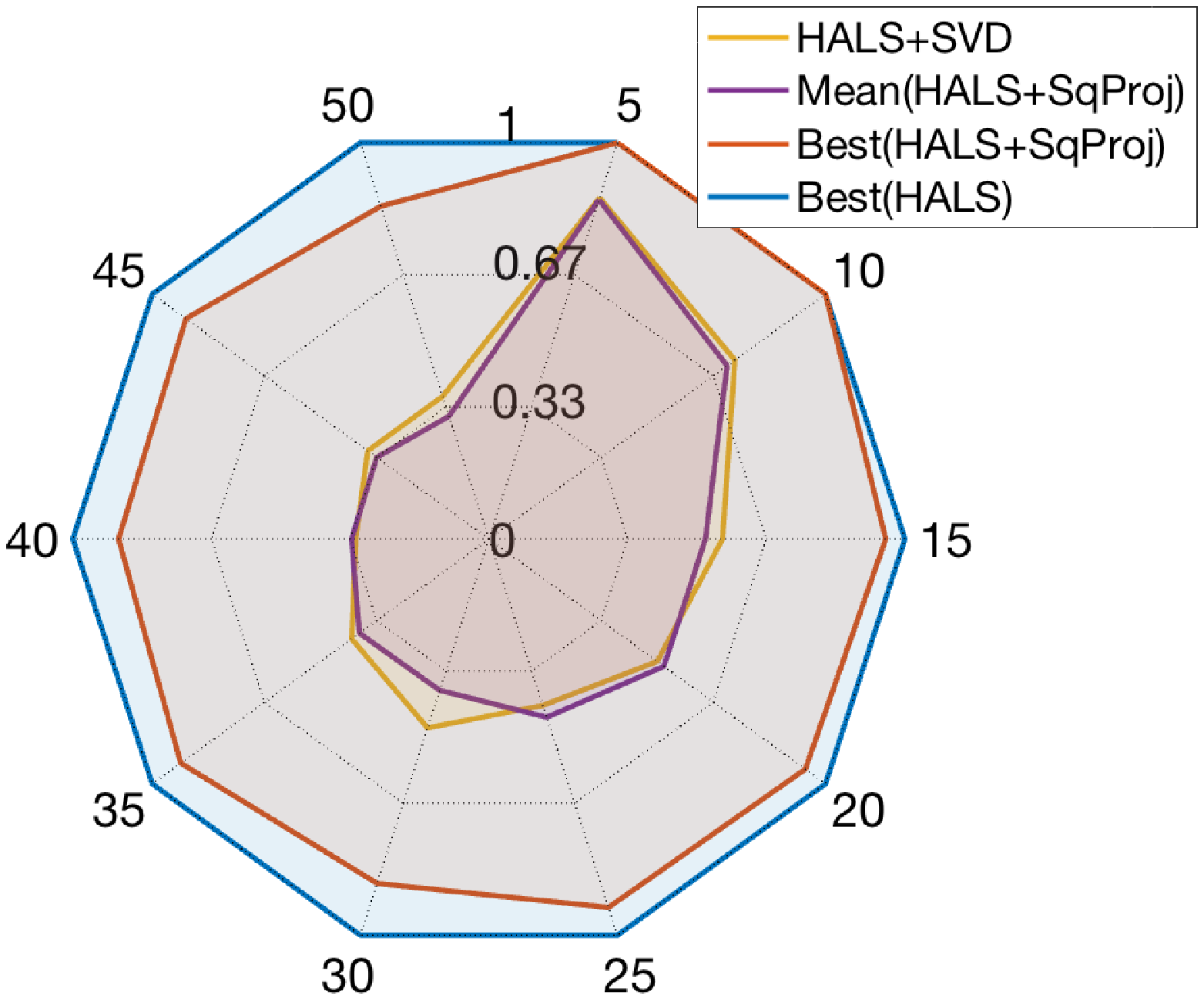}\label{fig_hals3}}
\subfigure[$N = 4$]{\includegraphics[width=.48\linewidth, trim = 1.0cm 1cm 0cm 0cm,clip=true]{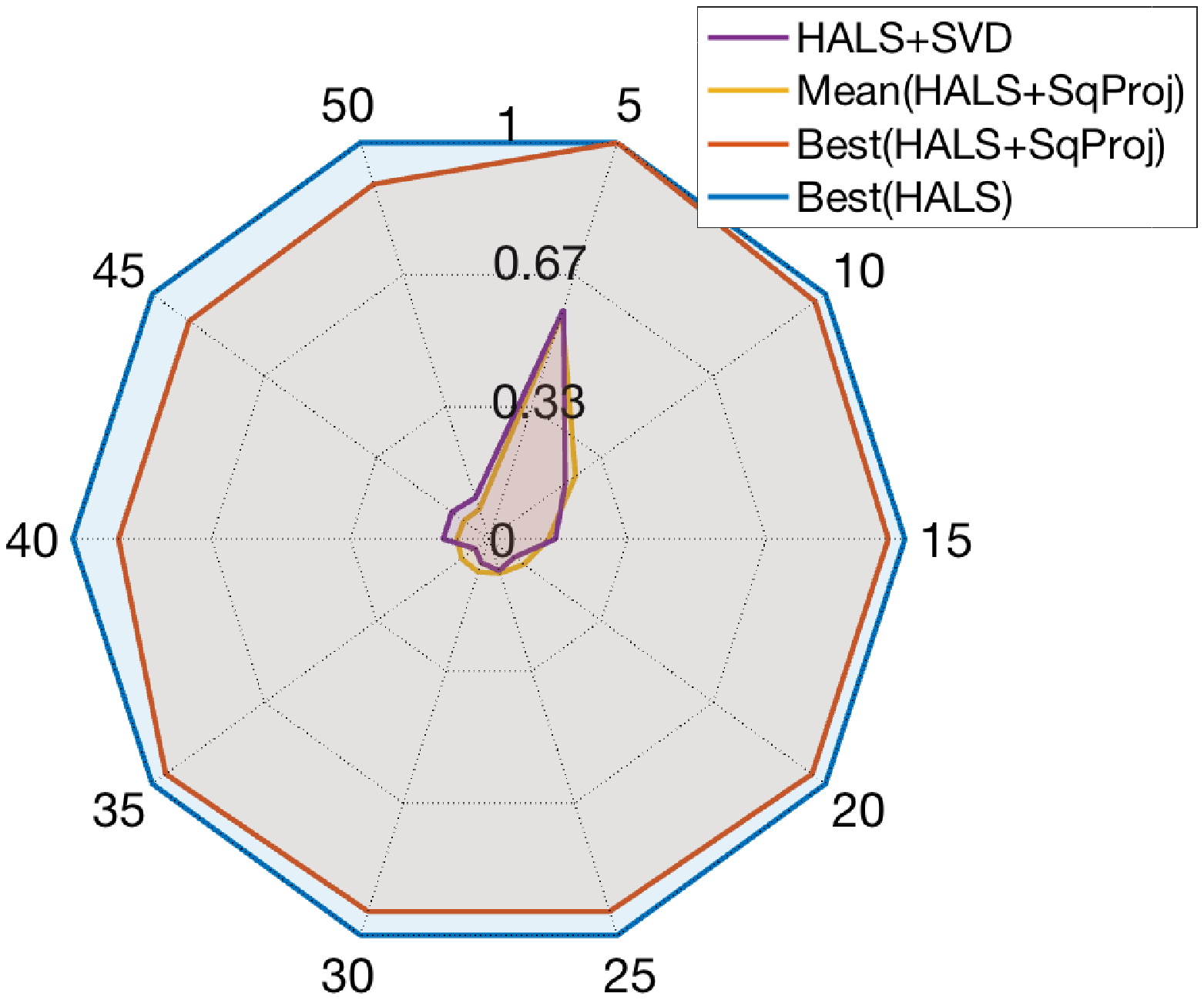}\label{fig_hals4}
}
\caption{Radar plots of the success ratios at $10^{-6}$ of algorithms in Example~\ref{ex_svd_tt_svd}, i.e., the percentage of approximation errors which were different from the best approximation errors with an error less than $10^{-6}$. 
Radius represents the success ratio, whereas angle corresponds to the tensor dimension $I = 5, 10, \ldots, 50$. A larger area indicates a higher success ratio.}
\label{fig_hals_comp_1e6}
\end{figure}

\begin{figure}[t!]
\centering
\subfigure[$N = 3$]{\includegraphics[width=.48\linewidth, trim = 1.0cm 1cm 0cm 0cm,clip=true]{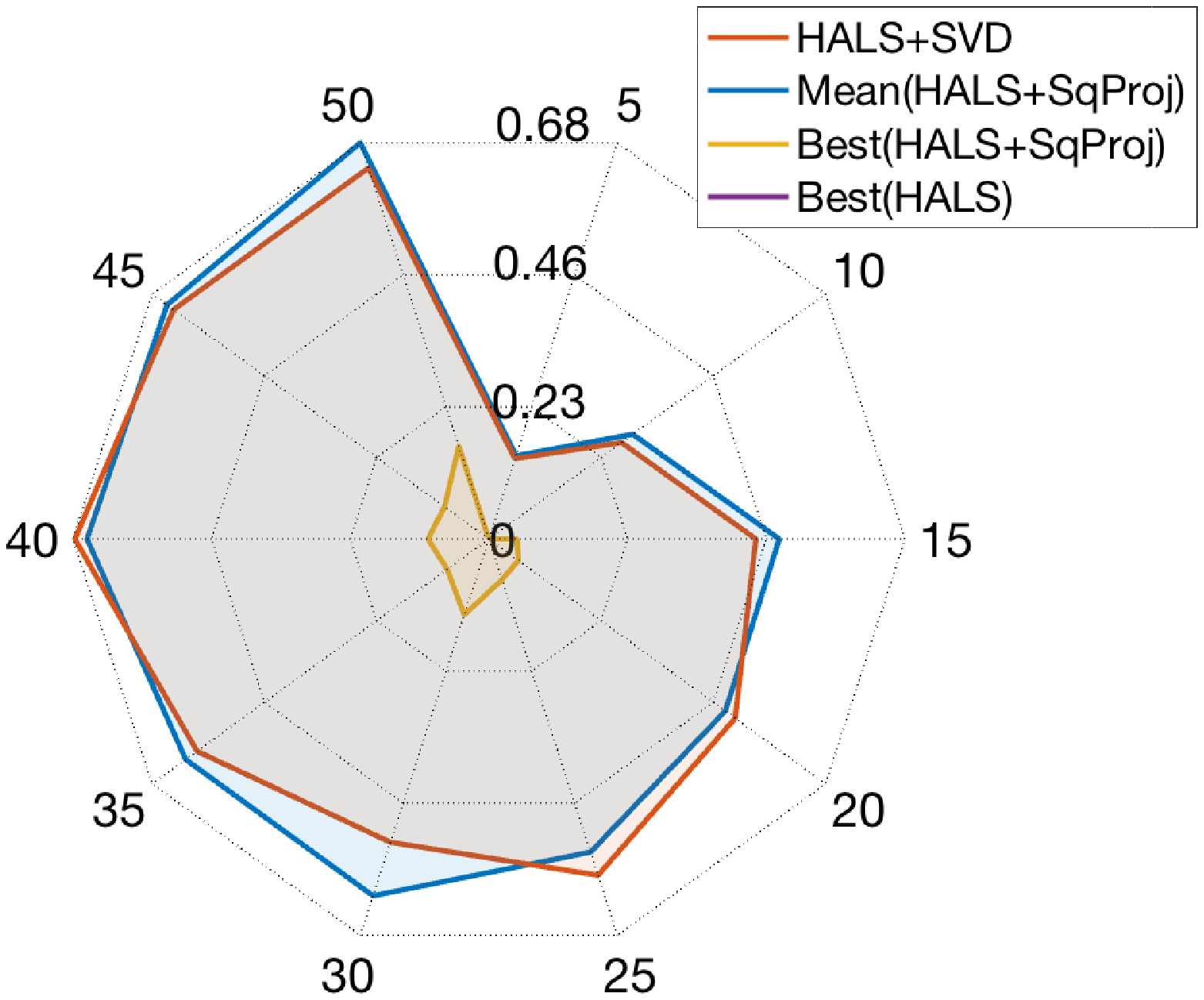}\label{fig_hals32}}
\subfigure[$N = 4$]{\includegraphics[width=.48\linewidth, trim = 1.0cm 1cm 0cm 0cm,clip=true]{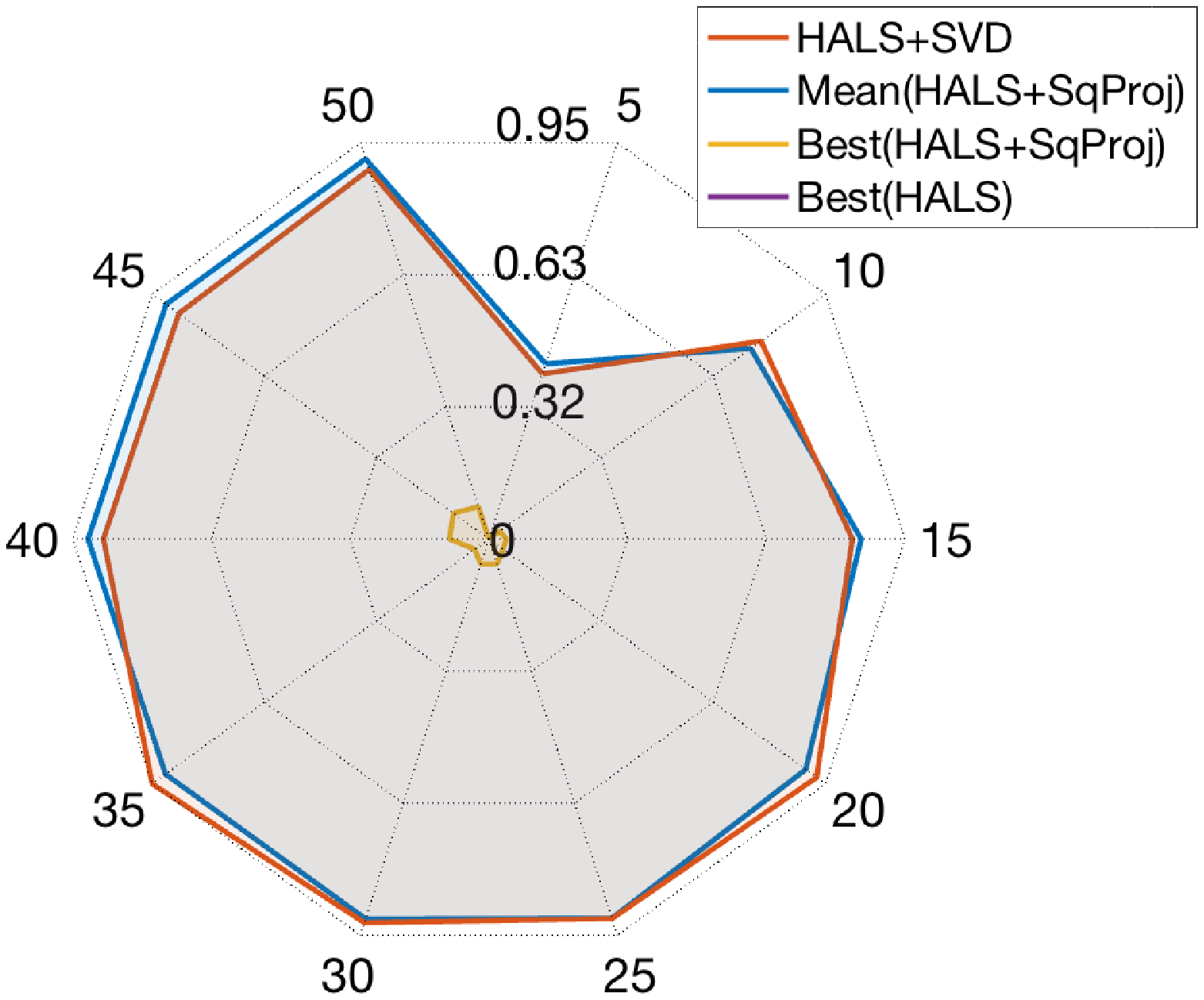}\label{fig_hals42}
}
\caption{Radar plots of the failure ratios at $10^{-2}$ of algorithms in Example~\ref{ex_svd_tt_svd}, i.e., the percentage of approximation errors which were different from the best approximation errors with an error greater than $10^{-2}$. 
Radius represents the failure ratio, whereas angle corresponds to the tensor dimension $I = 5, 10, \ldots, 50$.
A smaller area indicates a lower failure ratio.}
\label{fig_hals_comp_1e2}
\end{figure}

\begin{example}[SVD vs TT-SVD for best rank-1 tensor approximation]\label{ex_svd_tt_svd}

In this example, we seek the best rank-1 tensor approximation for random tensors of order-3 and 4, and tensor dimension $I = 5, 10, \ldots, 50$. The HALS was initialised using singular vectors, and TT-SVD applied to all $N!$ possible permutations of the tensor. The results were reported over at least 100 independent runs for each test case. 

For each run, the best approximation error was chosen among all $(N! + 1)$ results including one for HALS+SVD, and $N!$ errors for HALS using sequential projection (HALS+SqProj) for all $N!$ possible tensor permutations. 
We assessed percentage of approximation errors of an algorithm which were different from the best approximation error less than $10^{-6}$. This also represents the success ratio at $10^{-6}$.
These ratios are plotted in the radar plots in Fig.~\ref{fig_hals_comp_1e6}. 
\begin{itemize}

\item HALS using SVD initialisation achieved quite low success ratios, especially for tensors of order $N = 4$. 
The success ratios were lower for larger tensors.
Only for the case when tensors were of small size, $I = 5$, the success ratios of HALS+SVD were of 86.15\% and 58.00\% for $N = 3$ and $N = 4$, respectively.
\item The success ratios achieved by the sequential projection and truncation method (SqProj) were on average compatible with those using the SVD-based method. 
However, the best performances using the SqProj method, i.e., with proper tensor permutations, were much better than that of HALS+SVD. Its success ratios@$10^{-6}$ were respective of 92.38\% and 94.05\% for tensor orders $N = 3$ and $4$.
\end{itemize}

In Fig.~\ref{fig_hals_comp_1e2}, we illustrate the failure ratio at $10^{-2}$, i.e., the percentage of the approximation errors of an algorithm which were different from the best performance with an error greater than $10^{-2}$. On average the HALS+SVD failed to achieve the best results in 83.39\% and 50.02\% of runs for $N = 3, 4$.
\end{example}

Finally, despite that SVD-based and SqProj initialization methods are widely used for initialisation in tensor decompositions, the two methods often converge to local minima. 
With a suitable tensor permutation, SqProj may help to achieve the best approximation error.

 
\subsection{Levenberg-Marquardt algorithm for best rank-1 tensor approximation}

\subsubsection{Energy-balanced normalization} \label{sec:normalization}
Before introducing the proposed LM algorithm for best rank-1 tensor approximation, we present an energy-balanced normalisation for loading components.

Let $\gamma_n = \bu_n^T \bu_n$ and $\alpha = \sqrt[N]{ \gamma_1 \gamma_2 \cdots \gamma_N}$. 
Due to the scaling ambiguity, the loading components $\bu_n$ can be normalised to have balanced $\ell_2$-norm, i.e.,
\be
	\tilde{\bu}_n = \sqrt{\frac{\alpha}{\gamma_n}} \bu_n \label{eq_normalize_un}
\ee
to give $\tilde{\bu}_n^T \, \tilde{\bu}_n = \alpha$. 
This transformation preserves the rank-1 tensor 
\be
	\bu_1 \circ \bu_2 \circ \cdots \circ \bu_N = \tilde{\bu}_1 \circ \tilde{\bu}_2 \circ \cdots \circ \tilde{\bu}_N\,.
\ee

\subsubsection{The LM update}
We consider the following objective function
\be
\min  \quad \frac{1}{2} \|\tY - \bu_1 \circ \bu_2 \circ \cdots \circ \bu_N\|_F^2   + \frac{\mu}{2} \sum_{n = 1}^{N} \|\bu_n\|_2^2
\ee
 and apply the Levenberg-Marquardt algorithm to update the parameters $\btheta = [\bu_1^T, \ldots, \bu_N^T]^T$ 
\be
	\btheta \leftarrow \btheta -  (\bH + \mu \bI)^{-1} \, \bg \label{eq_LM}
\ee
where $\mu > 0$ is the damping parameter, $\bg$ and $\bH$ are gradient and approximate Hessian of the first term w.r.t. $\btheta$.
According to Theorem~2\cite{Phan_fLM}, the gradient and Hessian are given by   
\be
\bg &=&  \left[\ldots,  \bu_n  \gamma_{-n} - \bt_n ,\ldots \right]^T  \label{eq_grad_g}  \, ,\\ 
\bH_{\mu} &=& \bH + \mu \bI = \bD_{\mu} + \bZ \bK \bZ^T  \label{eq_H_1} \, ,
\ee
where   
\be
\bt_n &=& {\tY \, \bar{\times}_{k \neq n} \, \bu_k}, \label{eq_tn}
\\
\bD_{\mu} &=& \blkdiag( (\mu + \gamma_{-n})  \bI_{I_n}) \, ,\\
\bZ &=& \blkdiag(\ldots,   \bu_n, \ldots) \,, \\
\bK &=& [k_{n,m}] , \quad k_{n,n} = \0, \quad k_{n\neq m} = {\gamma_{-(n,m)}} 
\ee 
and $\gamma_{-n} = \prod_{k \neq n} \bu_k^T \bu_k$ and ${\gamma_{-(n,m)}} = \prod_{k \neq n,m} \bu_k^T \bu_k $.
Applying the energy-balanced normalisation in (\ref{eq_normalize_un}) to $\bu_n$ after each update, we have 
$\gamma_{-n}  = \alpha^{N-1}$ and ${\gamma_{-(n,m)}} = \alpha^{N-2}$.
The gradient and Hessian in (\ref{eq_H_1}) are then rewritten as 
\be
\bg &=&  \alpha^{N-1} \, \btheta - \bt \, \label{eq_g_2} \, ,\\
\bH_{\mu} &=& (\mu +  \alpha^{N-1})  \, \bI  + \alpha^{N-2}  \, \bZ \, (\1_N \1_N^T - \bI_{N}) \bZ^T  \notag \\
&=& (\mu +  \alpha^{N-1})  \, \bI   -  \alpha^{N-2} \, \bZ\,  \bZ^T  + \alpha^{N-2}  \, \btheta \btheta^T     \notag \\
&=& \bF_{\mu} + \alpha^{N-2} \, \btheta \, \btheta^T\,,\label{eq_H_F}
\ee
where $\bt = [\bt_1^T,\ldots, \bt_n^T]^T$ is a vector concatenated from $\bt_n$, and $\bF_{\mu} = \blkdiag((\mu +  \alpha^{N-1})  \, \bI  - \alpha^{N-2} \, \bu_n \bu_n^T)$ is a symmetric block diagonal matrix.
Inverse of $\bF_{\mu}$ is computed through inverses of its blocks as
\be
\bF_{\mu}^{-1} &=& \frac{1}{\mu +  \alpha^{N-1}} \blkdiag( (\bI  -   \frac{\alpha^{N-2}}{\mu +  \alpha^{N-1}} \, \bu_n \bu_n^T)^{-1}) \, \notag  \\
 &=& \frac{1}{\mu +  \alpha^{N-1}} \blkdiag( \bI_{I_n}  +   \frac{\alpha^{N-2}}{\mu} \, \bu_n \bu_n^T) \,  .
\ee
It is obvious to verify that 
\be
\bF_{\mu}^{-1} \, \btheta &=&  \frac{1}{\mu +  \alpha^{N-1}} \left[ \bu_{n}  +   \frac{\alpha^{N-1}}{\mu} \, \bu_n \right]_{n = 1}^{N} \,  
 = \frac{1}{\mu} \, {\btheta} \, , \label{eq_Ftheta}\\
\bF_{\mu}^{-1} \, \bt &=&  \frac{1}{\mu +  \alpha^{N-1}} \left[ \bt_{n}  +   \frac{\alpha^{N-2}}{\mu} \, \bu_n (\bu_n^T \bt_n) \right]_{n = 1}^{N} \,  
\notag \\
&=&  \frac{1}{\mu +  \alpha^{N-1}}  \, \bt + \frac{\alpha^{N-2} \xi}{\mu(\mu +  \alpha^{N-1})} \, \btheta \label{eq_Ft},
\ee
where $\xi = \langle \tY, \bu_1 \circ \bu_2 \circ \cdots \circ \bu_N \rangle = \bu_n^{T} \,  \bt_n$ is inner product of $\tY$ and its rank-1 approximation tensor.

From (\ref{eq_H_F}), inverse of the Hessian is given as a rank-1 update of the inverse of $\bF_{\mu}$ 
\be
\bH_{\mu}^{-1} &=&   \bF_{\mu}^{-1}  -  \frac{\alpha^{N-2}}{1  + \alpha^{N-2}\,  \btheta^T \,  \bF_{\mu}^{-1}  \, \btheta} \, \bF_{\mu}^{-1} \btheta \btheta^T  \bF_{\mu}^{-1} \notag\\
&=&  \bF_{\mu}^{-1} - \frac{\alpha^{N-2}}{\mu^2 + N \mu \alpha^{N-1}} \btheta \, \btheta^T \, . \label{eq_Hinv}
\ee

From (\ref{eq_g_2}), (\ref{eq_Ftheta}), (\ref{eq_Ft}) and (\ref{eq_Hinv}), the new increment vector $\bdelta = - \bH_{\mu}^{-1} \, \bg$ in the LM update (\ref{eq_LM}) is computed as
 \be
 \bdelta  &=&- (\bF_{\mu}^{-1} - \frac{\alpha^{N-2}}{\mu^2 + N \mu \alpha^{N-1}} \btheta \, \btheta^T) (\alpha^{N-1} \, \btheta - \bt)  \notag \\
 &=& \frac{-\alpha^{N-1} }{\mu} \btheta + \frac{N \alpha^{2(N-1)}}{\mu^2 + N \mu \alpha^{N-1} } \, \btheta  - \frac{N \alpha^{N-2} \xi}{\mu^2 + N \mu \alpha^{N-1}}  \btheta   \notag \\
 &&
 + \frac{\alpha^{N-2} \xi}{\mu(\mu +  \alpha^{N-1})} \, \btheta   +  \frac{1}{\mu +  \alpha^{N-1}}  \, \bt  \notag \\
 &=&  \frac{-\alpha^{N-1}}{\mu + N  \alpha^{N-1} } \, \btheta  - \frac{(N-1) \alpha^{N-2} \xi}{(\mu + N \alpha^{N-1})(\mu + \alpha^{N-1})}  \btheta    +  \frac{1}{\mu +  \alpha^{N-1}}  \, \bt  \notag  \\
 &=&  -\frac{\alpha^{N-1}(\mu +  \alpha^{N-1}) + (N-1) \alpha^{N-2} \xi}{(\mu + N \alpha^{N-1})(\mu + \alpha^{N-1})}  \btheta    +  \frac{1}{\mu +  \alpha^{N-1}}  \, \bt  \notag  \,.
 \ee
This finally leads to the update for $\btheta$ as follows
\be
\btheta \leftarrow  \frac{1}{\mu + \alpha^{N-1}} \left( \left(\mu + \frac{(N-1)\alpha^{N-2}  (\alpha^{N} - \xi)}{\mu+ N \alpha^{N-1}} \right)  \btheta + \bt \right) \,. \label{eq_LMupdate_theta}
\ee

\subsubsection{Simplication of the LM update with optimal norm of the rank-1 tensor}

Assume $\bu_1 \circ \bu_2 \circ \cdots \circ \bu_N$ be a rank-1 approximation tensor of the tensor $\tY$ after an LM update, 
we can always find a scaling factor $\lambda$ which minimises  the approximation error 
\be
	\min \quad \| \tY - \lambda \, \bu_1 \circ \bu_2 \circ \cdots \circ \bu_N \|_F^2 \,. \notag 
\ee
This gives  
\be	
\lambda = \langle \tY,  \bu_1 \circ \bu_2 \circ \cdots \circ \bu_N\rangle / \gamma  = \xi / \gamma \, .\label{eq_lambda} 
\ee
if $\lambda  = 1$, there no need an adjustment of the rank-1 tensor. Otherwise, if $\lambda \neq 1$, we obtain a new rank-1 tensor with a lower approximation error.
Adjusting the loading components $\bu_n$ by a factor of $\lambda^{1/N}$, i.e., 
\be
	\tilde{\bu}_n  = {\bu_n} \, \lambda^{1/N},  \notag 
\ee
 preserves the estimated rank-1 tensor
\be
	\lambda\,  \bu_1 \circ \bu_2 \circ \cdots \circ \bu_N =  \tilde{\bu}_1 \circ \tilde{\bu}_2 \circ \cdots \circ \tilde{\bu}_N \notag 
\ee
and the energy-balanced $\tilde{\bu}_n^T \tilde{\bu}_n = \lambda^{2/N} \alpha =\tilde{\alpha}$.
A more important result is that the new tensor has 
$\tilde{\xi} = \tilde{\gamma} = \tilde{\alpha}^{N}$ since 
\be
\tilde{\gamma} &=& \prod_n (\tilde{\bu}_n^T \, \tilde{\bu}_n) = \lambda^2 \gamma  \notag \,,\\
\tilde{\xi} &=& \langle \tY,  \tilde{\bu}_1 \circ \tilde{\bu}_2 \circ \cdots \circ \tilde{\bu}_N \rangle  =  \lambda \langle \tY,   {\bu}_1 \circ  {\bu}_2 \circ \cdots \circ  {\bu}_N \rangle  
=   \lambda^2 \, \gamma\, \notag. 
\ee
That is we can always convert a rank-1 tensor to have balanced energy and an optimal norm with $\xi = \gamma$.

Now we apply the LM update to such loading components with the equality $\xi = \gamma$, the LM update rule in (\ref{eq_LMupdate_theta}) becomes a simple update rule 
\be
\btheta &\leftarrow & \frac{1}{\mu + \alpha^{N-1}} \left( \mu   \btheta + \bt \right) \, \notag \\
&=&  \bu_n  -  \eta \, \bg_n  
\label{eq_update_un_from_unvn}  \label{eq_LMupdate_theta2}
\ee
where $\bg_n = \gamma_{-n}  \bu_n - \bt_n$  represent the gradient given in (\ref{eq_grad_g}) and 
\be
\eta &=& \frac{\displaystyle1}{\displaystyle \mu + \alpha^{N-1}}\ee
$\eta$ is considered a step size in the range of $\left[0, \frac{1}{\alpha^{N-1}}\right]$.

 \setlength{\algomargin}{1em}
\begin{algorithm}[t!]
\SetFillComment
\SetSideCommentRight
\CommentSty{\footnotesize}
\caption{{\tt{LM algorithm for best rank-1 tensor approximation}}\label{alg_lm_1}}
\DontPrintSemicolon \SetFillComment \SetSideCommentRight
\KwIn{Data tensor $\tY$:  $(I_1 \times I_2 \times \cdots \times I_N)$} 
\KwOut{$\tX =  \bu_1 \circ  \bu_2 \circ \cdots \circ \bu_N$}
\SetKwFunction{cpd}{CPD} 
\Begin{
\nl Initialize $\bu_n$ \;
\nl Normalise  $\bu_n \leftarrow  \frac{\alpha}{\sqrt{\gamma_n}} \,  \bu_n$, where $\alpha =   {\sqrt[N]{\xi}}:{\sqrt[2N]{\gamma}}$\;
\Repeat{a stopping criterion is met}{
\nl Solve the optimal step-size $\eta \in [0, \frac{1}{\alpha^{N-1}}]$ which minimises the polynomial $f(\eta)$ in (\ref{eq_f_eta})\;
\nl Update in parallel $\bu_n \leftarrow  \bu_n - \eta\, \bg_n$ \;
\nl Normalise $\bu_n \leftarrow  \frac{\alpha}{\sqrt{\gamma_n}} \,  \bu_n$, where $\alpha =   {\sqrt[N]{\xi}} \, : \sqrt[2N]{\gamma}$\;
}
}
\end{algorithm}

With the new form in (\ref{eq_update_un_from_unvn}), the LM  update rule reduces to the steepest-descend method, 
and finding the damping parameter $\mu$ is equivalent to seeking a step-size $\eta$. 
When the damping parameter is sufficient large, $\mu \rightarrow \infty$ ($\eta \rightarrow 0$), the above update rule cancels the term $\bt$, and there is no update here. In other words, we can choose a suitable $\mu$ to lower the approximation error. 
For example, the damping
parameter $\mu$ can be updated using Nielsen's method  \cite{NielsendampingLM}, or optimally determined as a root of a polynomial of degree-$2N$.

Note that due to the optimal norm condition, $\bg_n$ is orthogonal to $\bu_n$
\be
	\bg_n^T \bu_n &=& \alpha^{N-1} \bu_n^T \bu_n  -\bt_n^T \bu_n  
	 = \gamma - \xi  = 0 \, . \notag 
\ee
By exploiting this result,  we can rewrite the Frobenius norm of the error as a function of the step-size $\eta$
\be
f(\eta) &=& \|\tY - \hat{\tY}\|_F^2    \notag \\
&=& \|\tY\|_F^2  +  \|\hat{\tY}\|_F^2 - 2  \langle \tY, \hat{\tY}\rangle   \notag \\
&=&  \|\tY\|_F^2  +  \prod_{n = 1}^{N}  [1, \eta]^T [\bu_n, -\bg_n]^T  [\bu_n, -\bg_n]    [1, \eta]  
 - 2 \left\langle \tW ,  \left[\begin{array}{@{}c@{}}1\\[-1em]\eta\end{array}\right] \circ \left[\begin{array}{@{}c@{}}1\\[-1em]\eta\end{array}\right] \circ \cdots \circ \left[\begin{array}{@{}c@{}}1\\[-1em]\eta\end{array}\right] \right\rangle \notag \\
&=&  \|\tY\|_F^2  +  \prod_{n = 1}^{N}  (\alpha + c_n \, \eta^2)   - 2  \sum_{n = 0}^{N}  q_n \eta^n     \label{eq_f_eta}
\ee
where  $c_n = \bg_n^T\bg_n$ and $\tW = \tY \times_1 [\bu_1, -\bg_1]^T \times_2 [\bu_2, -\bg_2]^T  \cdots \times_N  [\bu_N, -\bg_N]^T$ is a tensor of size $2 \times 2 \times \cdots \times 2$.
 Inner product of the tensor $\tW$ and the vector $\left[\begin{array}{@{}c@{}}1\\[-1em]\eta\end{array}\right]$ yields a degree-$N$ polynomial with coefficients $q_n$.
  For simplicity we denote coefficients of the tensor $\tW$ by $\bw = [w_ 1, \ldots, w_ {2^N}] = \vtr{\tW}$. 
For $N= 3$
\be
 q_3 &=& w_8  , \quad q_2  =  w_4 + w_6 + w_7, \quad \,\notag \\ 
 q_0 &= &w_1, \quad q_1 = w_2 + w_3 + w_5  \notag \,.
\ee
For $N=4$,
\be
 q_4 &=& w_{16},  \quad q_3  =  w_8 + w_{12} + w_{14} + w_{15}, \notag \\
 \quad q_2  &=&  w_4 + w_6 + w_7 + w_{10} + w_{11} + w_{13}, \quad \,\notag \\ 
 q_0 &= &w_1, \quad q_1 = w_2 + w_3 + w_5 + w_9  \notag \,.
\ee
Coefficients $q_n$ for higher tensor order $N$ can be recursively deduced from those of lower order.

From the degree-$2N$ polynomial in (\ref{eq_f_eta}), the optimal step size $\eta$ minimises the objective function in $[0, \frac{1}{\alpha^{N-1}}]$.
This can be accomplished by finding roots of the derivative of the polynomial.	
A simple implementation of the LM algorithm is shown in Algorithm~\ref{alg_lm_1}.






\subsection{A Rotational algorithm}

Before deriving a new algorithm which updates loading components by rotations, we consider the derivation of the HOOI algorithm again, which minimises the objective function with loading components on spheres 
\be
\min \;\frac{1}{2} \|\tY - \beta\,  \bu_1 \circ \bu_2 \circ \cdots \circ \bu_N \|_F^2    \,, \;\;
\text{s.t.}\;\; \bu_{n}^T \, \bu_{n} = 1,  n = 1, \ldots, N\,. \notag 
\ee
This leads to $\beta = \tY \, \bar{\times}_{n = 1}^{N} \, \bu_n = \xi$, and the optimization problem becomes 
\be
	\text{max} \quad  \frac{1}{2} \, ( \tY \, \bar{\times}_{n = 1}^{N} \, \bu_n )^2  , \;\;  
	\text{s.t.}  \quad    \bu_{n}^T \, \bu_{n} = 1, \quad n = 1, \ldots, N . \label{equ_costYUt} 
\ee
For this new constraint optimisation, we construct a Lagrangian function 
\begin{align}	
 L(\bu_1, \ldots, \bu_N, \lambda_1, \ldots, \lambda_N) =   \frac{1}{2} \, ( \tY \, \bar{\times}_{n} \, \bu_n  )^2 -  \frac{1}{2}\sum_{n=1 }^{N}   {\lambda}_n \, \left(\bu_n^T \, \bu_{n} - 1 \right)    \notag 
\end{align}
where $\lambda_n$ are Lagrange multipliers.
The gradient of the Lagrangian with respect to $\bu_{n}$ is given by 
\be
	 \tilde{\bg}_n  =  \xi \, \bt_n \,   - \lambda_n \bu_n   \, ,\label{equ_gradLagrange}
\ee
where $\bt_n$ is defined in (\ref{eq_tn}).
Since $\bu_n^T \bu_n = 1$, by setting $\tilde{\bg}_n$ to zero, we obtain $\lambda_n  =  \bu_n^T \bt_n  \, \xi = \xi^2$, 
and 
\be
\tilde{\bg}_n = \xi \bt_n  - \xi^2 \bu_n\, .
\ee
Following the steepest descent method, the loading components $\bu_{n}$ can be updated as
\be
	\bu_n \leftarrow \bu_n + \eta \, \tilde{\bg}_n  \,, \label{equ_updateUn}
\ee
where the step size $\eta >0$. 
The above update rule, however, may not preserve the unit-length constraints of $\bu_n$. 
Note that the new estimate of $\bu_n$ lies in the subspace spanned by $[\bu_n, \tilde{\bg}_n]$. Let $\bar{\bg}_n$ be unit-length vector of $\tilde{\bg}_n$, i.e., $\tilde{\bg}_n = \|\tilde{\bg}_n\| \, \bar{\bg}_n$. 
We rewrite the update rule in (\ref{equ_updateUn}) in a new form as
\be
	\bu_n \leftarrow [\bu_n, \bar{\bg}_n] \,  \boldeta_n \,, \label{equ_updateUn2}
\ee
where $\boldeta_n = [\eta_{n1}, \eta_{n2}]^T$ are rotational vectors of length 2.
Similar to $\bg_n$ in (\ref{eq_update_un_from_unvn}), it can be verified that 
\be
	\bu_n^T \, \tilde{\bg}_n =  (\bu_n^T \bt_n) \xi  - (\bu_n^T \bu_n) \xi^2 = 0.
\ee
Hence $[\bu_n, \bar{\bg}_n]^T [\bu_n, \bar{\bg}_n] = \bI_2$.
In order to preserve the unit-length constraints of $\bu_n$, the vectors $\boldeta_n$ must lie on a unit sphere, i.e., $\boldeta_n^T \boldeta_n = 1$ for all $n$.
We replace $\bu_n$ in (\ref{equ_costYUt}) by their new updates in (\ref{equ_updateUn2}), and find rotational vectors $\boldeta_n$ in a best rank-1 tensor approximation 
\be
	&&\text{max} \quad  \frac{1}{2} ( \tW \, \bar{\times}_{n = 1}^{N} \, \boldeta_n )^2 \label{equ_cost_rot_eta} \\
	&&\text{s.t.}  \quad    \boldeta_n^T \, \boldeta_n = 1, \quad n = 1, \ldots, N  \notag  
\ee
where
$\tW$ is a projected tensor of size $2 \times 2 \times \cdots \times 2$ from the tensor $\tY$
\be
\tW = \tY\, \bar{\times}_1 [\bu_1, \bar{\bg}_1]  \, \bar{\times}_2 [\bu_2, \bar{\bg}_2]   \cdots \bar{\times}_N  [\bu_N, \bar{\bg}_N] \,  \label{eq_W_proj} \,.
\ee 
A simple implementation of the above update rules is listed in Algorithm~\ref{alg_r1rotational}.
In each update,  the ROtational algorithm to find best Rank-One tensor  (RORO) algorithm seeks best rank-1 tensor approximation to quantised-scale tensors of size $2 \times 2 \times \cdots \times 2$.
which can solved using the ALS/HOOI or R1LM algorithm.
We will show that this step can be done in closed-form for tensors of size $2 \times 2 \times 2$, and best rank-1 tensor of higher order tensors can be found efficiently from those of lower order.

 \setlength{\algomargin}{1em}
\begin{algorithm}[t!]
\SetFillComment
\SetSideCommentRight
\CommentSty{\footnotesize}
\caption{\tt{ROtational Algorithm for best Rank-One tensor approximation (RORO)}\label{alg_r1rotational}}
\DontPrintSemicolon \SetFillComment \SetSideCommentRight
\KwIn{Data tensor $\tY$:  $(I_1 \times I_2 \times \cdots \times I_N)$} 
\KwOut{$\tX =  \bu_1 \circ  \bu_2 \circ \cdots \circ \bu_N$}
\SetKwFunction{cpd}{CPD} 
\Begin{
\nl Initialize unit length vectors $\bu_n$ \;
\Repeat{a stopping criterion is met}{
\nl Compute the projected tensor $\tW = \tY \times_1 [\bu_1, \bar{\bg}_1]^T \times_2 [\bu_2, \bar{\bg}_2]^T  \cdots \times_N  [\bu_N, \bar{\bg}_N]^T$\;
\nl Seek unit-length vectors $\boldeta_1$, \ldots, $\boldeta_N$ in best rank-1 tensor approximation to $\tW$\; 
\nl Update in parallel all $\bu_n \leftarrow  [\bu_n, \bar{\bg}_n]  \boldeta_n$ \;
}
}
\end{algorithm}

\subsubsection{RORO and best rank-1 tensor approximation to tensors of size $2 \times 2 \times 2$}

When the tensor $\tY$ is of order-3, the projected tensor $\tW$ is of size $2 \times 2 \times 2$. 
It is obvious that $\boldeta_1$ and $\boldeta_2$ are leading singular vectors of the projected matrix, $\tW \bar{\times}_3 \boldeta_3$, of size $2 \times 2$. 
From the problem in (\ref{equ_cost_rot_eta}), maximising the product $(\tW \bar{\times}_{n=1}^{3} \boldeta_n)^2$ is equivalent to maximizing the largest singular value of the matrix $\tW \bar{\times}_3 \boldeta_3$, i.e., 
\be
	\max  \quad  \sigma_{\max}^2(\tW \bar{\times}_3 \, \boldeta_3)  \,. \label{eq_obj_eta3}
\ee

 By repamaterizing $\boldeta_3 = \left[\cos(\alpha), \sin(\alpha)\right]^T$ and representing the $2 \times 2$  projected martrix $\bW(:,:,1) \cos(\alpha) + \bW(:,:,2) \sin(\alpha)$, 
following Appendix~\ref{sec_svd_2x2}, we can formulate the objective function in (\ref{eq_obj_eta3}) as a new optimisation to find $\alpha$ in $[0, 2\pi]$
\be
\max \quad f(\alpha) = a(\alpha) + \sqrt{a^2(\alpha) - b^2(\alpha)}\,  \label{eq_opt_rot_eta}\, ,
\ee
where 
 \be
 a(\alpha) &=&  a_1  \cos(2\alpha)  + a_2   \sin(2\alpha) + a_3 \, ,\\
 b(\alpha) &=&  b_1  \cos(2\alpha)  + b_2   \sin(2\alpha) + b_3   \, ,
 \ee
where the parameters, $a_k$  and $b_k$, are provided in Appendix~\ref{sec::ak_bk}.

By changing the parameter $\alpha = \arctan(x)$, and after some manipulations, the maximiser $\alpha^{\star}$ to $f(\alpha)$ in (\ref{eq_opt_rot_eta}) can be found by solving a degree-6 polynomial equation
\be
p(x) =  c_6 x^6 + c_5 x^5 +  \cdots  + c_1 x + c_0  = 0\label{eq_poly_6}
\ee
where the coefficients $c_k$ are given in Appendix~\ref{sec:poly6}. 
Among real-valued roots, we choose the root $x^{\star}$ associated with the largest value $f(\arctan(x^{\star}))$ in (\ref{eq_opt_rot_eta}).

 \setlength{\algomargin}{1em}
\begin{algorithm}[t!]
\SetFillComment
\SetSideCommentRight
\CommentSty{\footnotesize}
\caption{\tt{Closed-form for Best Rank-1 Tensor of a $2 \times 2 \times 2$ Tensor}\label{alg_222}}
\DontPrintSemicolon \SetFillComment \SetSideCommentRight
\KwIn{Data tensor $\tY$:  $(2 \times 2 \times 2)$} 
\KwOut{$\tX =  \sigma \, \bu_1 \circ  \bu_2 \circ  \bu_3 \approx \tY$}
\SetKwFunction{cpd}{CPD} 
\Begin{
\nl  Find roots of a degree-6 polynomial $p(x) = c_6 x^6 + \ldots + c_1 x + c_0$ with coefficients defined in (\ref{eq_c6})-(\ref{eq_c0})\; 
\nl  Choose the root $x^{\star}$ which yields largest $f(\arctan(x^{\star}))$ in (\ref{eq_opt_rot_eta})\;
}
\nl $\bu_3 = \frac{1}{\sqrt{1+x^2}} [1, x]^T$ \;
\nl Find a best rank-1 $\tY \bar{\times}_3 \bu_3 \approx \sigma \, \bu_1 \bu_2^T$\;
\end{algorithm}

 In summary, we compare the ALS/HOOI and the two new algorithms R1LM and RORO in Table~\ref{tab_compare_alsroro}.
 
\begin{table*}[t!]
\centering
\caption{Comparison of algorithms for best rank-1 tensor approximation to tensors of size $I \times I \times I$.}
\label{tab_compare_alsroro}
\begin{tabular}{lccc}
& ALS/HOOI & R1LM & RORO \\\hline
Update & One component per iteration &  All 3 components&  All 3 components in closed-form  \\
Subproblem & none & \multicolumn{2}{c}{Root of a polynomial of degree-6}  \\
Optimal $\lambda$ in (\ref{eq_lambda}) & yes & correct $\lambda$ after updating $\bu_n$  & yes \\
Computational cost & $3I^3$ & $5 I^3$  & $5 I^3$ \\
	& (due to $\bt_n$) & (due to $\bt_n$ and $\tW$) & (due to $\bt_n$ and $\tW$) \\
Initialisation & \multicolumn{3}{c}{Sequential projection with best tensor permutation} \\ \hline
\end{tabular}
\end{table*}

  \begin{figure}[t!]
\centering
\includegraphics[width=.45\linewidth, trim = 0.0cm 0cm 0cm  0cm,clip=true]{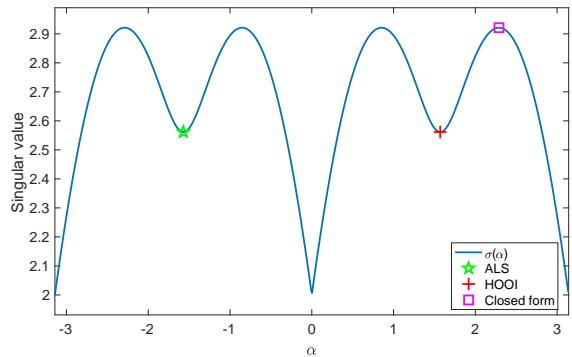}\label{fig_ex_222}
\caption{Singular values of the projected matrix $\tY  \bar{\times}_3 \boldeta_3$ in Example~\ref{ex_222}. ALS/HOOI fails to attain the largest singular value which is associated with the best rank-1 tensor.}
\label{fig_bestrank222}
\end{figure}

 \begin{example}[Failure of HOOI and ALS.]\label{ex_222}
We decompose a simple $2 \times 2 \times 2$ tensor $\tY$ whose entries are given by 
 \be
 \bY_1 = \left[\begin{array}{rr} 
   0 &   2\\
   2 &  0
 \end{array}\right]\, , \quad 
  \bY_2 = \left[\begin{array}{rr} 
     0 &   2\\
   -2 &  -1
 \end{array}\right] \, .
 \ee
Fig.~\ref{fig_bestrank222} shows the largest singular values, $\sigma(\alpha)$ in (\ref{eq_opt_rot_eta}), of the projected matrix  $\tY \bar{\times}_3 \, [\cos(\alpha), \sin(\alpha)]^T$ for $\alpha$  varied in the interval $[-\pi, \pi]$.
The best rank-1 tensor of $\tY$ has a scaling factor of $\tY \, \bar{\times}_{n = 1}^{3} \, \boldeta_n = 2.9212$.
Using the SVD-based initialisation method, the ALS/HOOI algorithm converges to a rank-1 tensor with a scaling factor of 2.5616. 
Besides the above tensor, ALS/HOOI often fails in seeking the best rank-1 tensor of the following tensors whose two frontal slices, [$\bY_1$, $\bY_2$], are given by
\be
&&\left[\left[\begin{array}{cc} 
   2 & -2 \\
   1 & 0
 \end{array}\right]\,,\left[\begin{array}{cc} 
   0 &   0\\
   2   &  2
 \end{array}\right] \right] \, ,
 \left[\left[\begin{array}{cc} 
   1 & -1 \\
   2 & -1
 \end{array}\right]\,,\left[\begin{array}{cc} 
   1  &    -2\\
  -2  &  -2
 \end{array}\right] \right]\, ,
%
%
 \left[\left[\begin{array}{cc} 
   -2 & 1 \\  
   1 & -2
 \end{array}\right]\,,\left[\begin{array}{cc} 
   -1  &  2\\
   0  &  2
 \end{array}\right] \right]\, \, .
%
\ee  \notag 
 \end{example}
            
\subsubsection{RORO and best rank-1 tensor approximation to tensors of size $2 \times 2 \times 2 \times 2$}
 
We can solve the problem in (\ref{equ_cost_rot_eta}) for tensors of order-4 following the alternating update scheme. At each iteration, we compress the tensor $\tW$ by a vector $\boldeta_n$ to yield an order-3 tensor $\tW \bar{\times}_n \, \boldeta_n$.
Then, the three unit-length vectors $\boldeta_k$, $k \neq n$ are found in closed-form as best rank-1 tensor of this projected tensor. The algorithm proceeds the estimation of another three components from the other compressed tensors until a convergence is achieved. The procedure is listed in Algorithm~\ref{alg_2222}.
Practical simulation results indicate that the algorithm executes only a few iterations.
The method can be straightforwardly extended to higher order tensors.

 \setlength{\algomargin}{1em}
\begin{algorithm}[t!]
\SetFillComment
\SetSideCommentRight
\CommentSty{\footnotesize}
\caption{\tt{Best Rank-1 Tensor of a $2 \times 2 \times 2 \times 2$ Tensor}\label{alg_2222}}
\DontPrintSemicolon \SetFillComment \SetSideCommentRight
\KwIn{Data tensor $\tY$:  $(2 \times 2 \times 2  \times 2)$} 
\KwOut{$\tX =  \bu_1 \circ  \bu_2 \circ \bu_3 \circ \bu_4 \approx \tY$}
\SetKwFunction{cpd}{CPD} 
\Begin{
\nl Initialize a unit length vector $\bu_1$ \;
\Repeat{a stopping criterion is met}{
\For{$n = 1, 2, 3, 4$}{
\nl  Apply Algorithm~\ref{alg_222} to find best rank-1 tensor of $2 \times 2 \times 2$ tensor $\tY \bar{\times}_n \bu_n$  \;
\nl Update component $\bu_k$, $k \neq n$\;
}
}
}
\end{algorithm} 

Alternatively, we can show that the parameters can be found as roots of a bi-variate polynomial of degree-6.
Let unit length vectors $\boldeta_3 = [\cos(\alpha), \sin(\alpha)]^T$ and  $\boldeta_4 = [\cos(\beta), \sin(\beta)]^T$, then from (\ref{eq_largest_sv2x2}), the largest singular value
of the projected matrix $\tW \,\bar{\times}_3 \, \boldeta_3 \,\bar{\times}_4 \, \boldeta_4 $ is  given by 
\be
	f(\alpha,\beta) = 2  \, \sigma_{\max}^2 = {a(\alpha,\beta) + \sqrt{a^2(\alpha,\beta) - b^2(\alpha,\beta)}} \label{eq_largest_sv2x2_b}
\ee
where 
\be
a(\alpha,\beta) &=&   [\cos(2\alpha), \sin(2\alpha), 1]   \, \bA \, [\cos(2\beta), \sin(2\beta), 1]^T \, ,  \notag \\
b(\alpha,\beta) &=&   [\cos(2\alpha), \sin(2\alpha), 1]   \, \bB \,  [\cos(2\beta), \sin(2\beta), 1]^T\,,  \notag 
\ee
the two matrices, $\bA$ and $\bB$, are of size $3 \times 3$ and defined in Appendix~\ref{sec:AB_2x2x2x2}.

Similarly to solving (\ref{eq_opt_rot_eta}), we perform a reparameterization $\alpha = \arctan(x)$ and $\beta = \arctan(z)$. 
The maximiser $(\alpha^{\star}, \beta^{\star})$ is root of the following bi-variate polynomial equations which maximises $	f(\alpha,\beta)$  
\be
p_1(x,z) =    [x^{6}, x^{5}, \ldots, x, 1] \, \bC_1  \,  [z^{6}, z^{5}, \ldots, z, 1]^T = 0\,,  \label{eq_bp_41_a}\\
p_2(x,z)  =    [x^{6}, x^{5}, \ldots, x, 1] \, \bC_2  \,  [z^{6}, z^{5}, \ldots, z, 1]^T = 0 \,, \label{eq_bp_42_a} 
\ee
where $\bC_1$ and $\bC_2$ are two matrices of size $7 \times 7$ derived in Appendix~\ref{sec:bivariate_poly}.

We can estimate $x$ as a root of a polynomial of degree-6, $p(x) = [x^{6}, x^{5}, \ldots, x, 1] \bC_1 [z^{6}, z^{5}, \ldots, z, 1]^T$,  while keeping $z$ fixed, then estimate $z$ in a similar way as a root of the polynomial  $q(z) = [x^{6}, x^{5}, \ldots, x, 1] \bC_2 [z^{6}, z^{5}, \ldots, z, 1]^T$ while $x$ is fixed.  
This is similar to the sequential projection of the tensor $\tW$ by either the vector $\boldeta_4$ or $\boldeta_3$. However, coefficients of the polynomials of degree-6 are simply provided either by $\bC_1 [z^{6}, z^{5}, \ldots, z, 1]^T$ or $[x^{6}, x^{5}, \ldots, x, 1] \bC_2$.

\begin{example}[Best rank-1 tensor approximation to a $2 \times 2 \times 2 \times 2$ tensor]\label{ex_bestrank1_2x2x2x2}
We decompose a simple $2 \times 2 \times  2 \times 2$ tensor $\tY$ whose four frontal slices are given by 
 \be
 \bY_{1,1} &=& \left[\begin{array}{rr} 
-1  &  -2\\
-2  &  2
 \end{array}\right]\, , \quad 
  \bY_{2,1} = \left[\begin{array}{rr} 
   0  &  -2 \\
   2  &   0
 \end{array}\right] \,,  \notag \\
   \bY_{1,2} &=& \left[\begin{array}{rr} 
    0  &   2 \\
-2     & 1
 \end{array}\right] \,,\quad
    \bY_{2,2} = \left[\begin{array}{rr} 
 -1   &  0 \\
 -2   & 1
  \end{array}\right] \,. \notag 
 \ee
Fig.~\ref{fig_bestrank2222} illustrates the largest singular value of the $2\times 2$ matrix projected from the tensor $\tY$ by the two unit length vectors $\boldeta_3 = [\cos(\alpha), \sin(\alpha)]^T$ and $\boldeta_4 = [\cos(\beta), \sin(\beta)]^T$.
Similar to Example~\ref{ex_222}, the ALS/HOOI algorithm fails to retrieve the best rank-1 tensor of the tensor $\tY$. ALS3 sequentially estimates three components at a time. In addition, we provide result obtained by solving the bi-variate polynomials in (\ref{eq_bp_41}) and (\ref{eq_bp_42}) using the IRIT multivariate solver\footnote{\url{http://www.cs.technion.ac.il/~gershon/irit/matlab/}}. ALS3 and IRIT achieve the best result,  
but they approach different points $(\alpha, \beta)$. In spite of that, the vectors $\boldeta_3$ and $\boldeta_4$ obtained by the two methods are the same after a sign correction. 
 
  \begin{figure}[t!]
\centering
\includegraphics[width=.45\linewidth, trim = 0.0cm 0cm 0cm  0cm,clip=true]{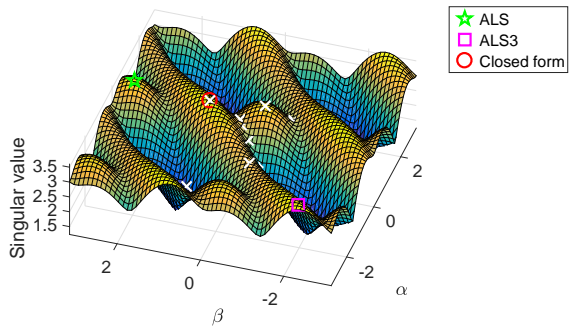}\label{fig_ex_2222}
\hspace{1em}
\caption{Singular values of the projected matrix $\tY  \bar{\times}_3 \boldeta_3 \bar{\times}_4 \boldeta_4$ in Example~\ref{ex_bestrank1_2x2x2x2}. ALS which updates components one-by-one converges to a false best rank-1 tensor, ALS3 which updates 3 components at a time yields the best rank-1 tensor. Roots of 
the bi-variate polynomials in (\ref{eq_bp_41}) and (\ref{eq_bp_42}) are denoted by ``$\times$'' markers. 	 
}  
\label{fig_bestrank2222}
\end{figure}
\end{example}

\section{Numerical Results}

\begin{example}[Best rank-1 tensor approximation]\label{ex_bestrank1}

We seek best rank-1 tensor approximations to random tensors of  order $N$ = 3, 4  and size $I \times I \times \cdots \times I$, where $I = 5, 10, \ldots, 30$. 
Four algorithms including HOOI, ALS, R1LM and RORO were initialised using leading singular vectors of each mode. Algorithms stopped when the changes in the approximation errors were lower than $10^{-12}$ or the number of iterations exceeded 1000.

\begin{figure}[t!]
\centering
\subfigure[Success ratio at $10^{-6}$.]{\includegraphics[width=.45\linewidth, trim = 0.0cm 0cm  0cm 0cm,clip=true]{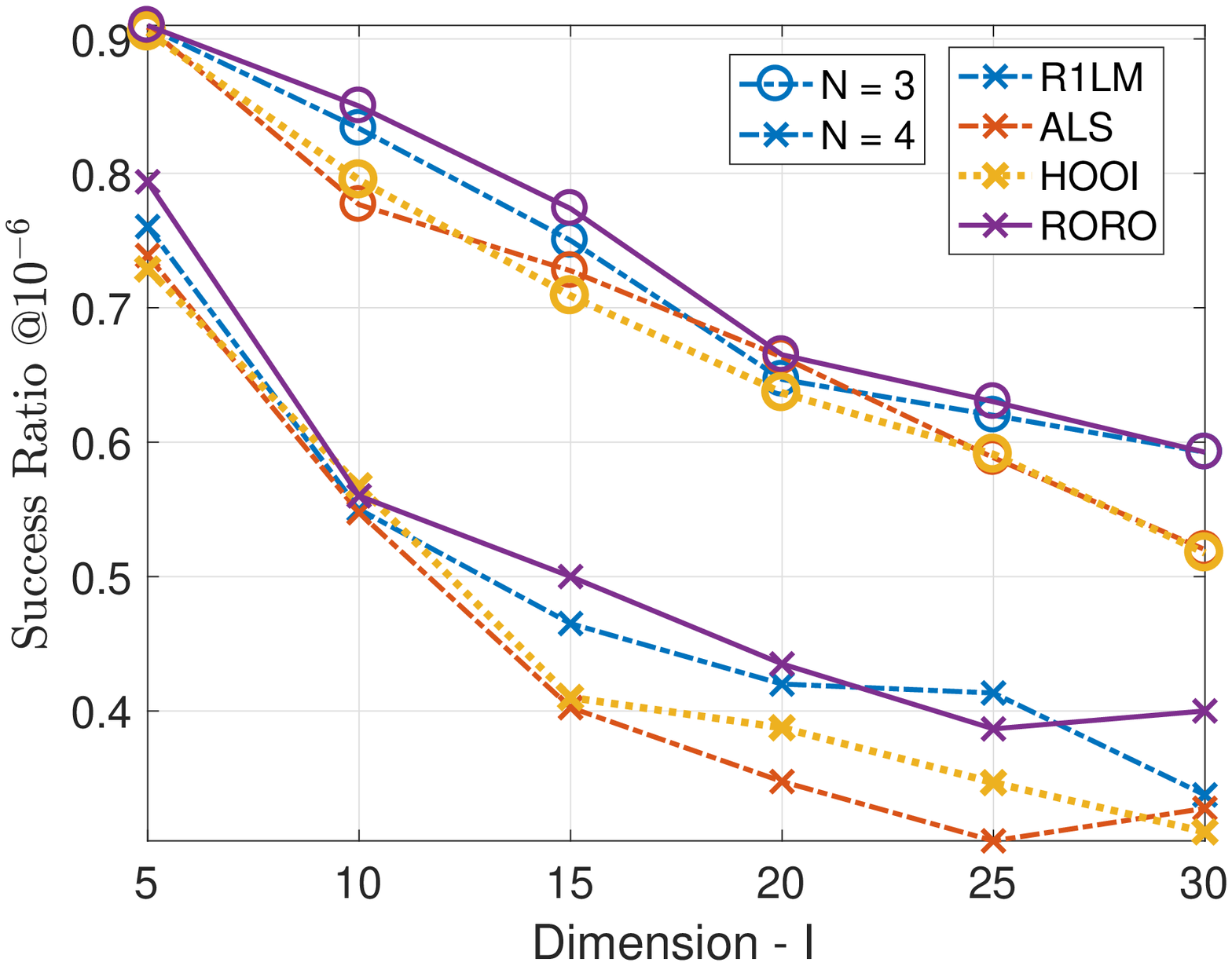}\label{fig_ex_t5d_a}}
\subfigure[Failure ratio at $10^{-3}$.]{\includegraphics[height=.35\linewidth, trim = 0.0cm 0cm 0cm 0cm,clip=true]{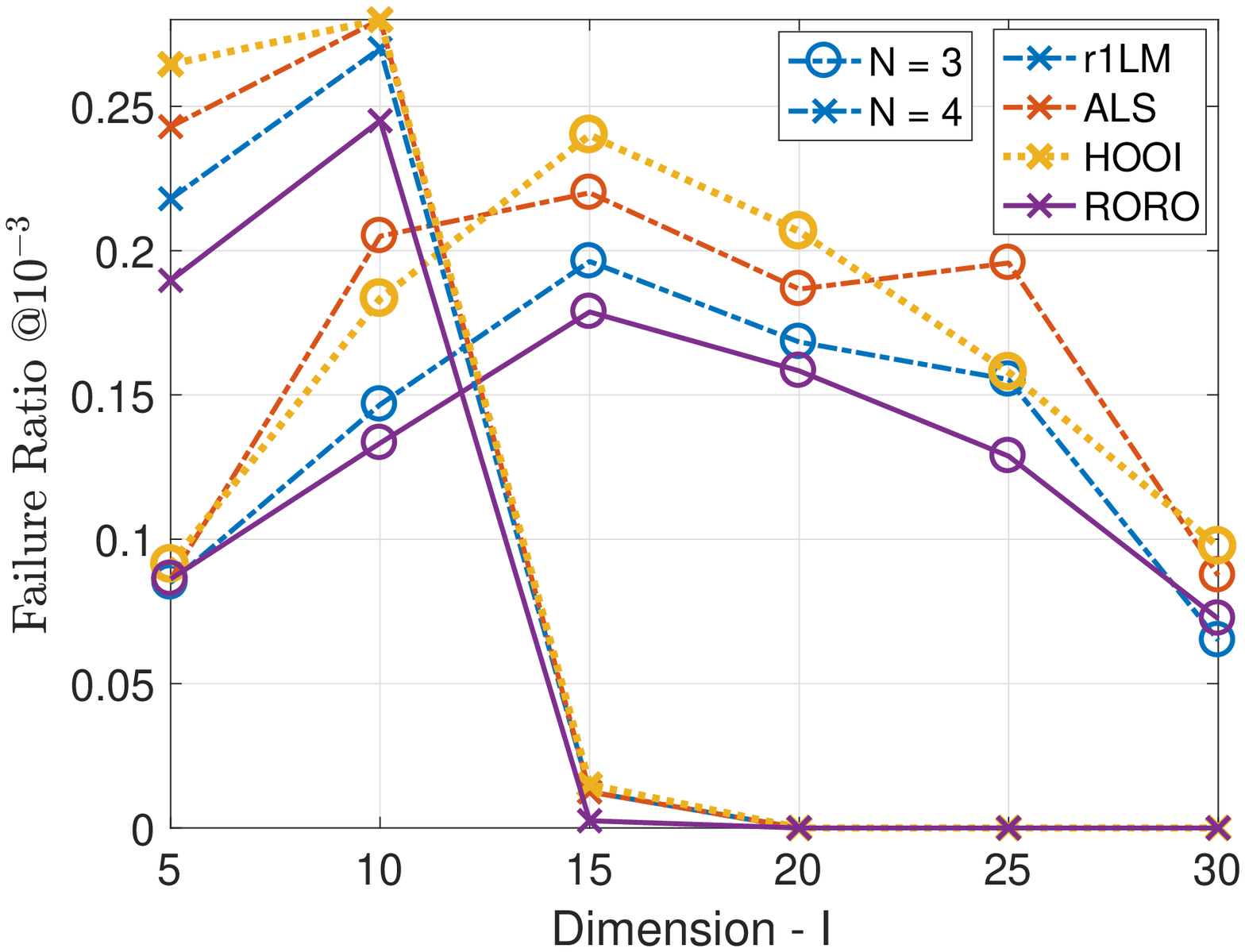}\label{fig_ex_t5d_b}
}
\caption{Comparison of success and failure ratios of algorithms for best rank-1 tensor approximation of order-3 and 4. Markers indicate the tensor order, $N$, whereas curves of an algorithm are plotted with the same line style. Algorithms were initialized using the SVD-based method.}
\label{fig_bestrank1_t5d_a}
\end{figure}

\begin{figure}[t!]
\centering
\subfigure[Success ratio at $10^{-4}$.]{\includegraphics[width=.45\linewidth, trim = 0.0cm 0cm 0cm 0cm,clip=true]{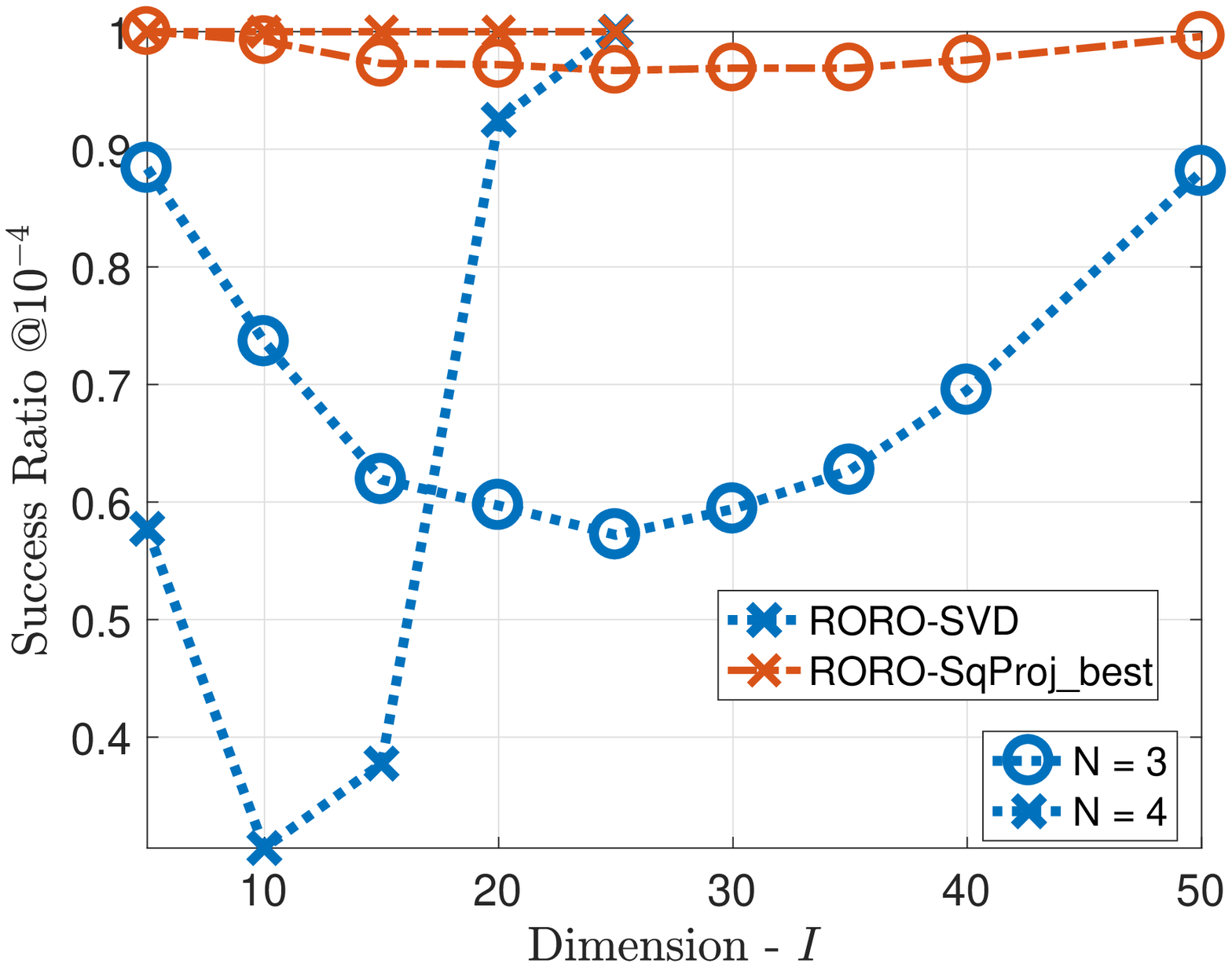}\label{fig_ex_t6_b}}
\subfigure[Success ratio at $10^{-6}$.]{\includegraphics[width=.45\linewidth, trim = 0.0cm 0cm 0cm 0cm,clip=true]{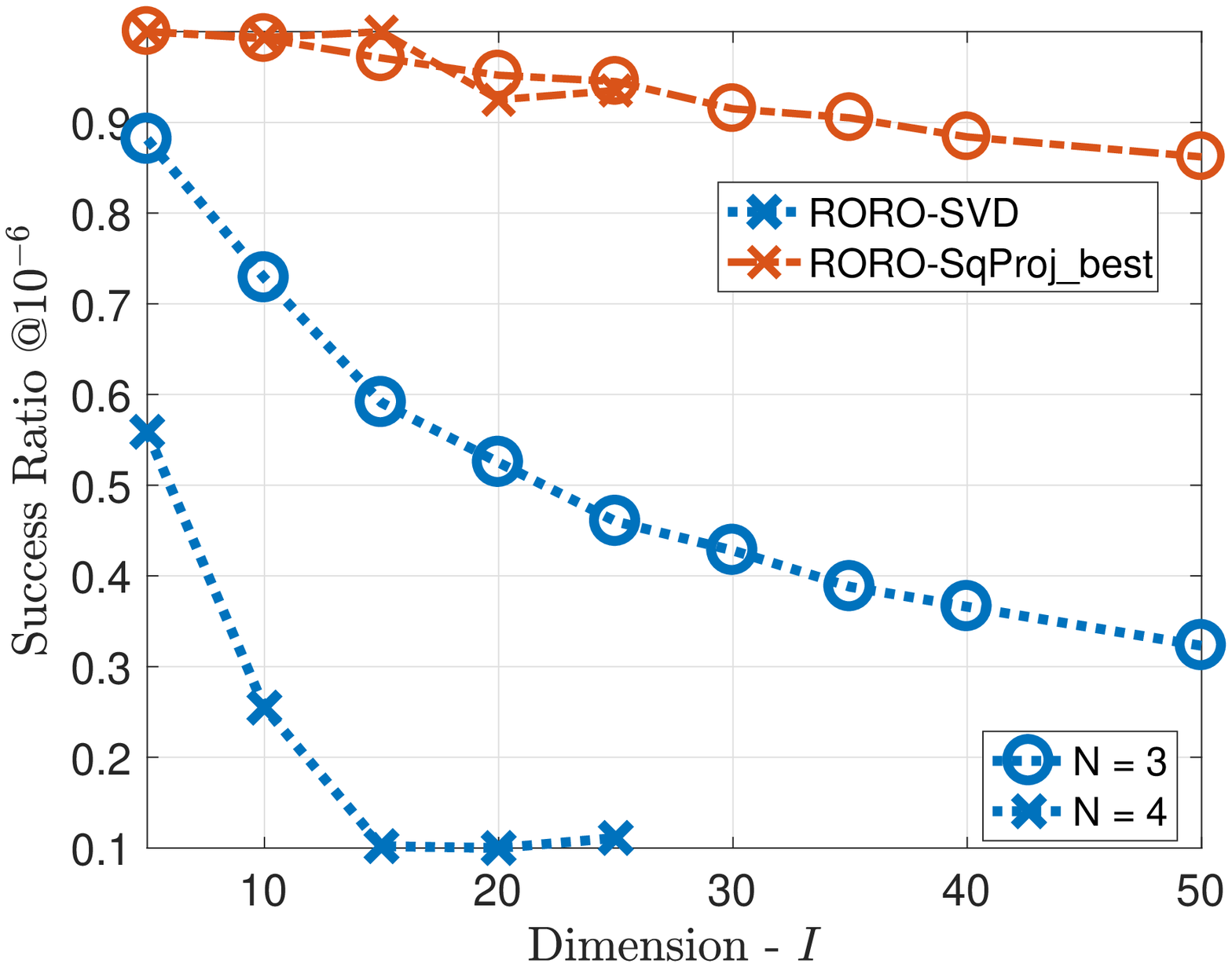}\label{fig_ex_t6_a}}
\caption{Success ratios of RORO using SVD and SqProj initializations. Markers indicate the tensor order, $N$, whereas curves of an algorithm are plotted with the same line style.}
\label{fig_bestrank1_t6}
\end{figure}

For each pair of $(N,I)$, we decomposed at least  400 tensors. Success ratios at $10^{-6}$ and failure ratios at $10^{-3}$ are compared in Fig.~\ref{fig_bestrank1_t5d_a}.

\begin{itemize}
\item In theory, ALS and HOOI have the same update rules. However, since ALS implements the fast projection, the update order of the loading components in the two algorithms are different \cite{Phan_fastALS}, their performances were slightly different. 
\item RORO and R1LM often achieved higher success ratios than the ALS/HOOI algorithms.
However, these two algorithms did not achieve perfect success ratios,  except for the case of small tensors of size $5 \times 5 \times 5$. A major reason is that the SVD-based initialisation method might lead to local minima. This can be improved e.g., using the sequential projection and truncation method as seen in Fig.~\ref{fig_bestrank1_t6}.
\item Algorithms converged to false local minima in at most 8 to 25 \% of runs for tensors of order-3, while R1LM and RORO have smaller failure ratios.
For tensors of order-4 and size $I \ge 15$, the failure ratios at $10^{-3}$ were almost zeros, although the estimated tensors were often be different from the best one with an error greater than $10^{-6}$.
\end{itemize}

In another comparison, success ratios at $10^{-4}$ and $10^{-6}$ of RORO using the sequential projection and truncation are plotted in Fig.~\ref{fig_bestrank1_t6}. The results indicates that the RORO+SqProj were much better than RORO using the SVD-based initialisation method. The results of R1LM were compatible to those of RORO.

\end{example}

\begin{example}[Decomposition of the multiplication tensors $(2 \times 2) \times (2 \times 2)$ and $(2 \times 3) \times (3 \times 2)$]
\label{ex_multiplicationtensor}

\begin{figure}[t!]
\centering
\subfigure[PARO with a fixed $\gamma$.]{\includegraphics[width=.45\linewidth, trim = 0.0cm  0cm 0cm .65cm,clip=true]{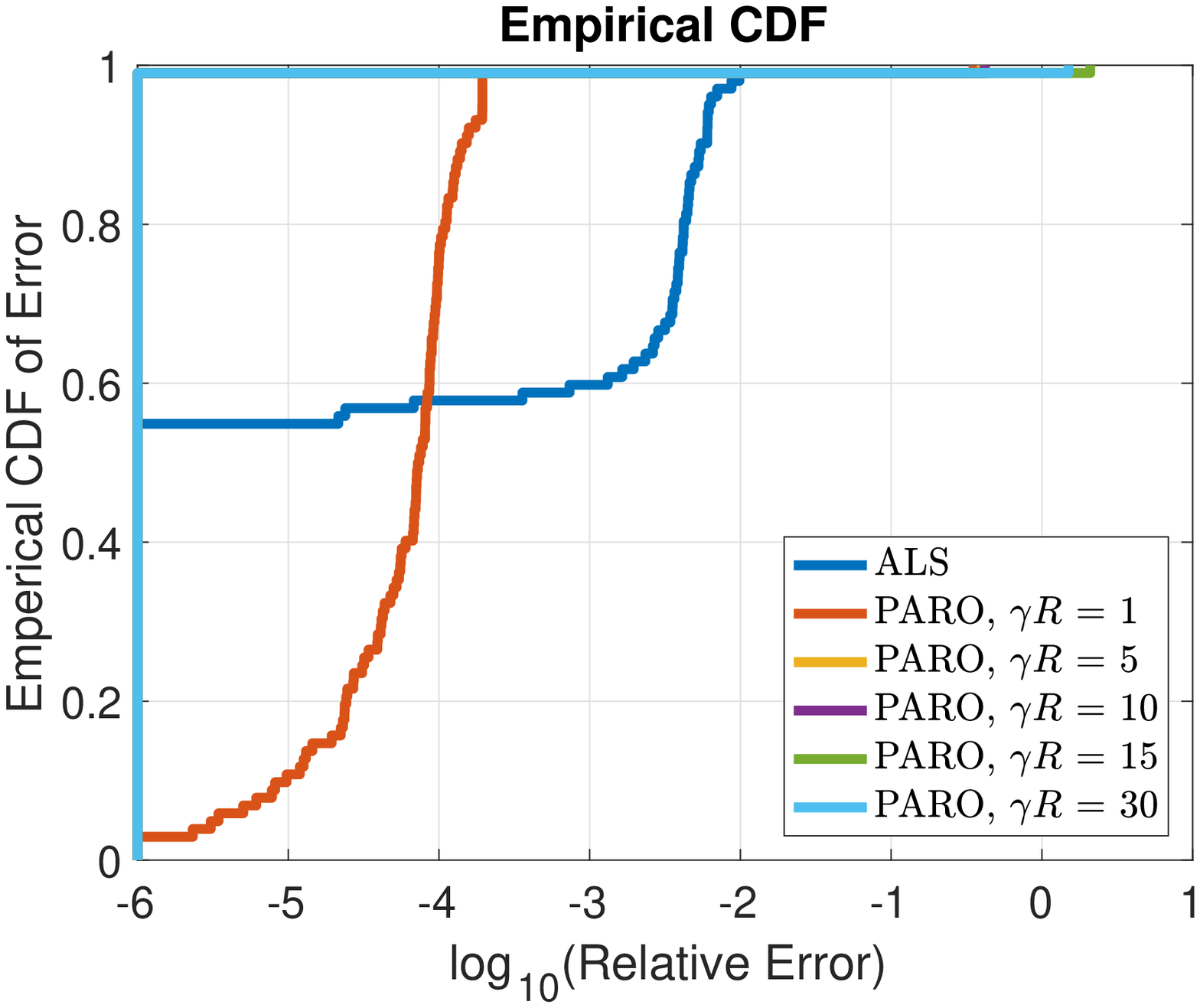}\label{fig_222b}
}
\subfigure[PARO with a regularly adjusted $\gamma$.]{\includegraphics[width=.45\linewidth, trim = 0.0cm 0cm 0cm .65cm,clip=true]{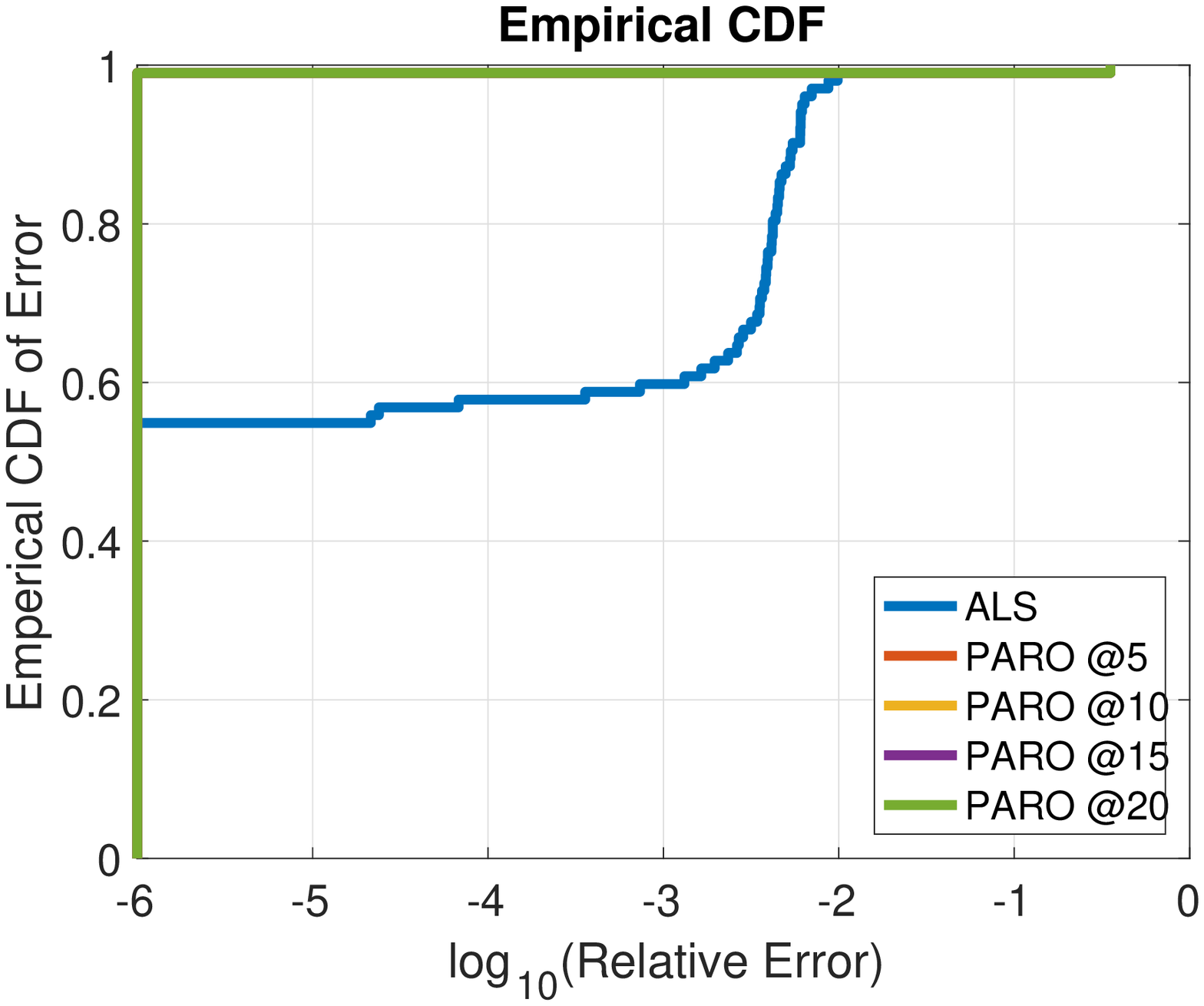}\label{fig_222a}}
\caption{Performance of PARO in decomposition of tensors of rank-$R = 7$ associated with the multiplication of two matrices $2\times 2$ and $2 \times 2$. Curves of PARO with $\gamma R = 5, 10, 15, 30$, or with $\gamma$ regularly adjusted are overlapped.}
\label{fig_2x2x2}
\end{figure}

\begin{figure}[t!]
\centering
\subfigure[PARO with fixed $\gamma$.]{\includegraphics[width=.45\linewidth, trim = 0.0cm  0cm 0cm .65cm,clip=true]{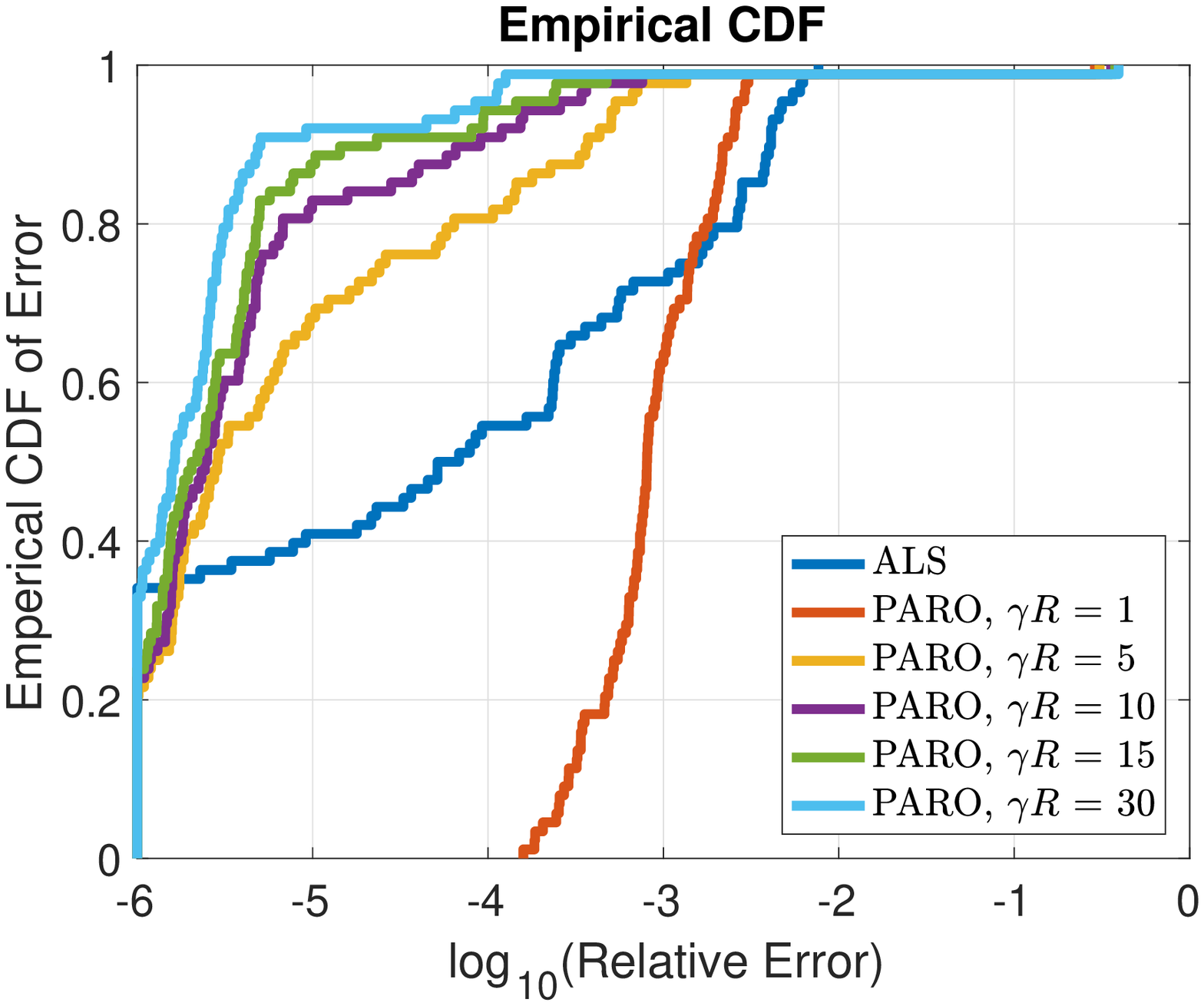}\label{fig_232_b}
}
\subfigure[PARO with regularly adjusted $\gamma$.]{\includegraphics[width=.45\linewidth, trim = 0.0cm 0cm 0cm .65cm,clip=true]{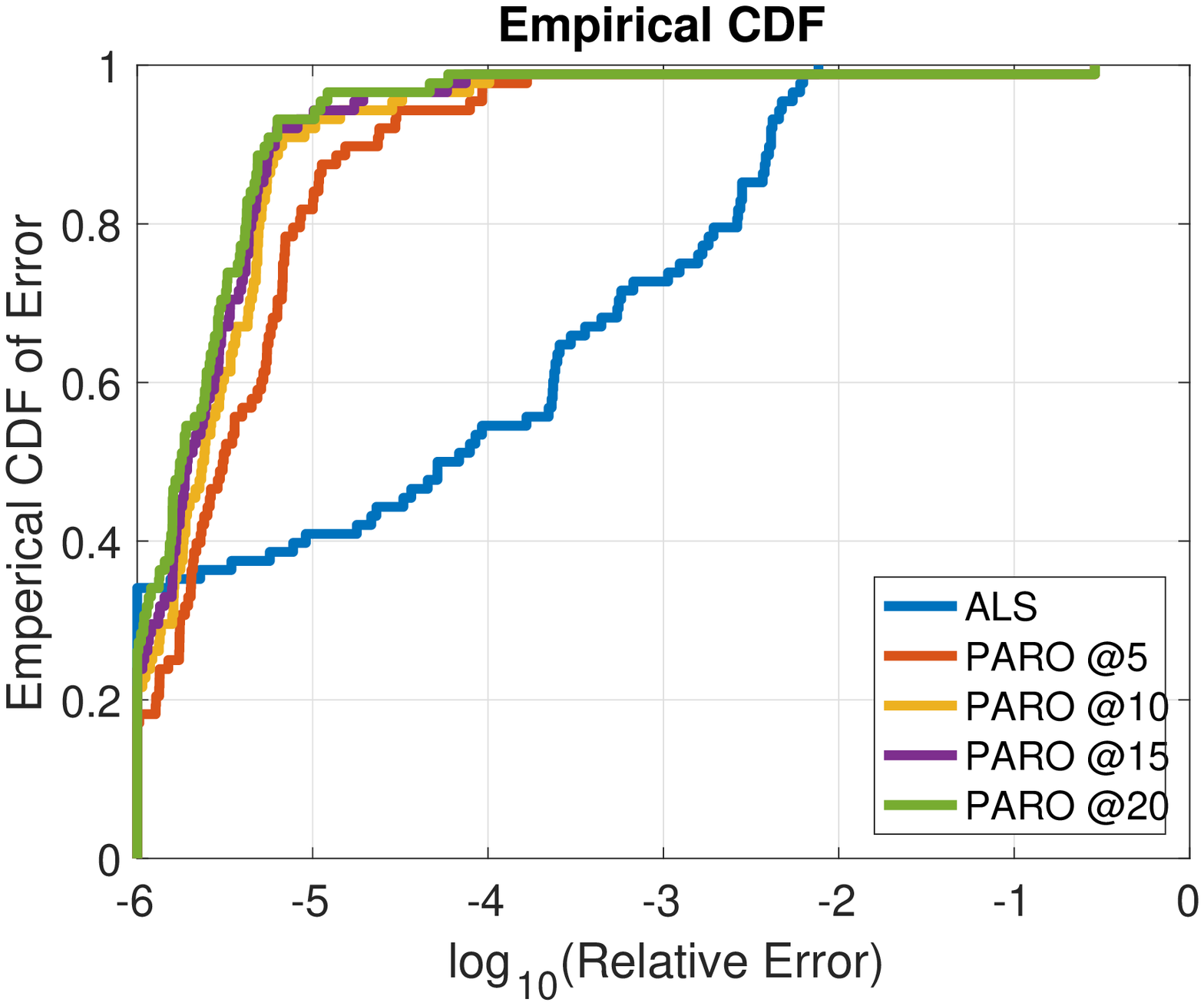}\label{fig_232a}}
\caption{Performance of PARO in decomposition of tensors of rank-$R = 7$ associated with the multiplication of two matrices $2\times 3$ and $3 \times 2$ in Example~\ref{ex_multiplicationtensor}.}
\label{fig_2x3x2}
\end{figure}

In this example, we decomposed multiplication tensors as in Example~\ref{ex_nasob222}. The tensors are associated with multiplication of two matrices of size $2 \times 2$, $2\times 2$, and $2 \times 3$, $3 \times 2$. The tensor rank of the latter tensor is $R = 11$.
We ran PARO using RORO with a fixed regularization $\gamma$ in 1000 iterations and with $\gamma$ adjusted every 5, 10, 15 and 20 iterations.
ALS was run in 5000 iterations but might stop if its relative error was smaller than $10^{-10}$.
Performances of the two algorithms are compared and shown in Fig.~\ref{fig_2x2x2} and Fig.~\ref{fig_2x3x2}.

\begin{itemize}
\item ALS could explain the tensors but in less than 60\% of runs for the tensors of size $4 \times 4 \times 4$, 
and even less than 40\% of runs for the tensors of size $4 \times 9 \times 4$.
\item PARO with the default regularization parameter $\gamma = 1/R$ had lower success ratios because the algorithm did not converge in 1000 iterations, but it will converge with more than 3000 iterations, as seen in Fig.~\ref{fig_222_1}.
\item PARO with higher regularisation parameters  $\gamma$ or with an adaptively adjusted gamma achieved almost perfect results.
\end{itemize}

\end{example}

\begin{example}[Decomposition of the tensor for the multiplication $(3 \times 3) \times (3 \times 3)$]

\begin{figure}[ht!]
\centering
\subfigure[PARO with fixed $\gamma$.]{\includegraphics[width=.45\linewidth, trim = 0.0cm  0cm 0cm .65cm,clip=true]{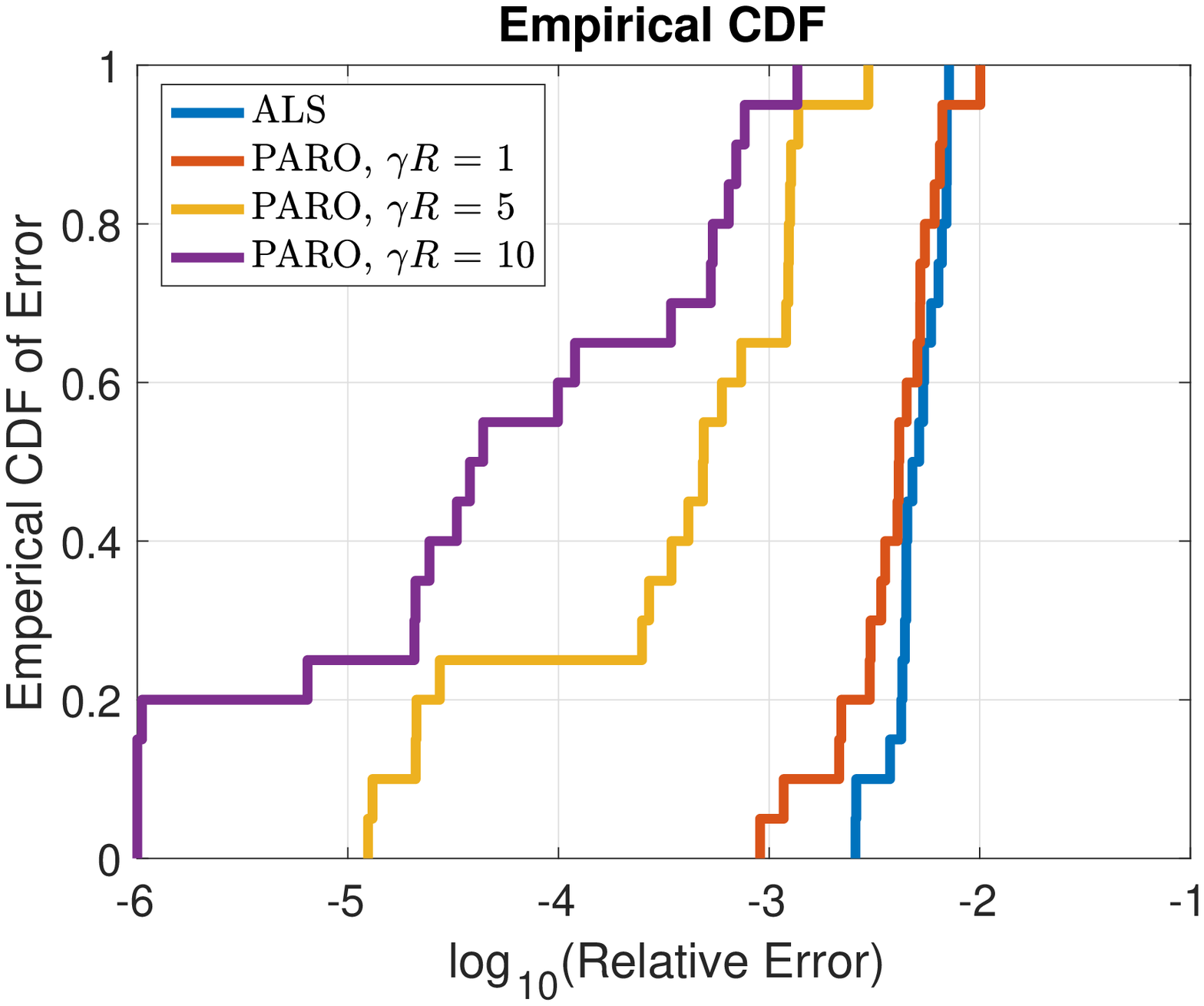}\label{fig_333a}
}
\subfigure[PARO with regularly adjusted $\gamma$.]{\includegraphics[width=.45\linewidth, trim = 0.0cm 0cm 0cm .65cm,clip=true]{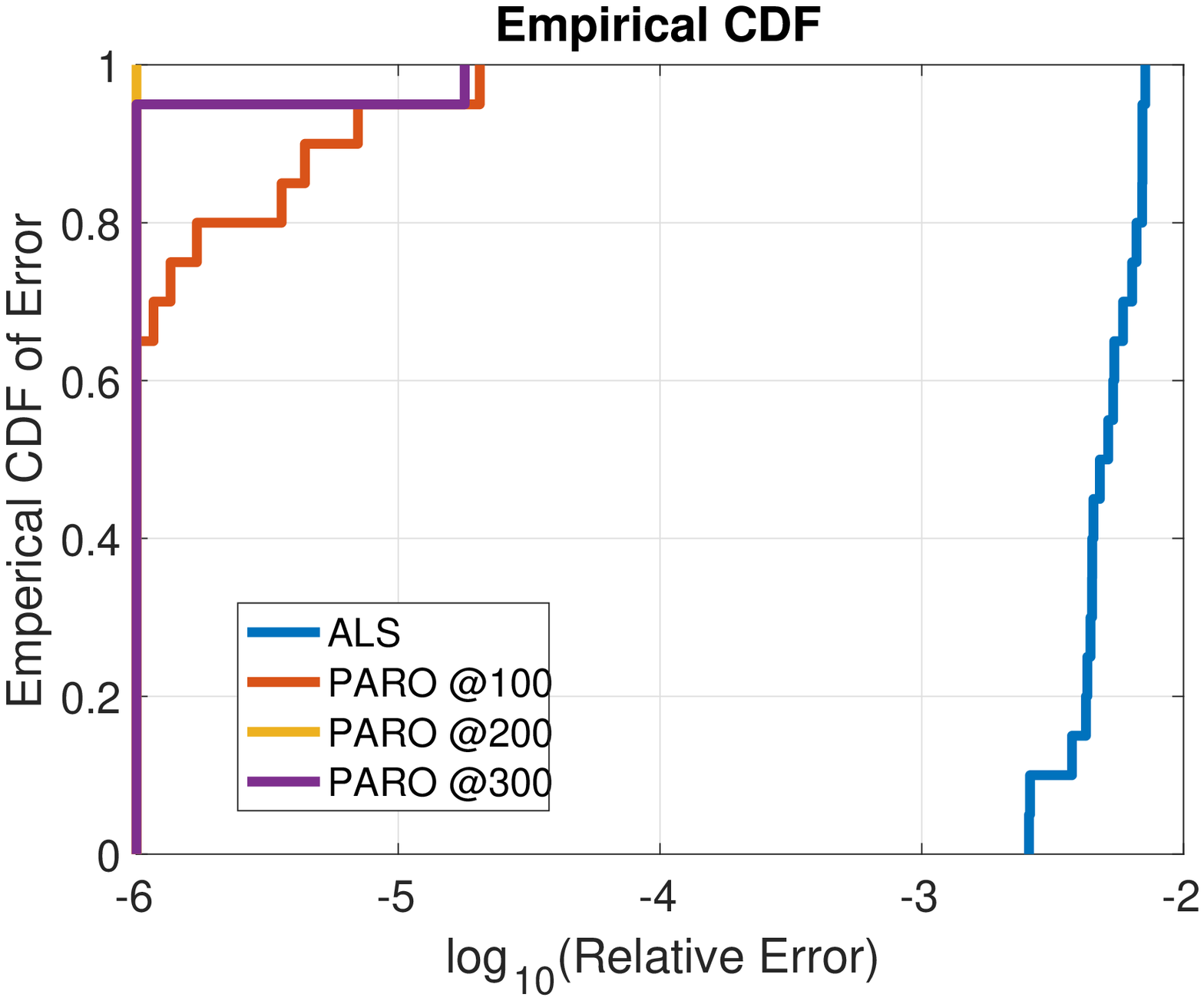}\label{fig_333b}}
\caption{Performance of PARO in decomposition of tensors of rank-$R = 23$ associated with the multiplication of two matrices $3\times 3$ and $3 \times 3$.}
\label{fig_3x3x3}
\end{figure}

\begin{figure}[ht!]
\centering
\subfigure[PARO with fixed $\gamma$.]{\includegraphics[width=.45\linewidth, trim = 0.0cm  0cm 0cm .0cm,clip=true]{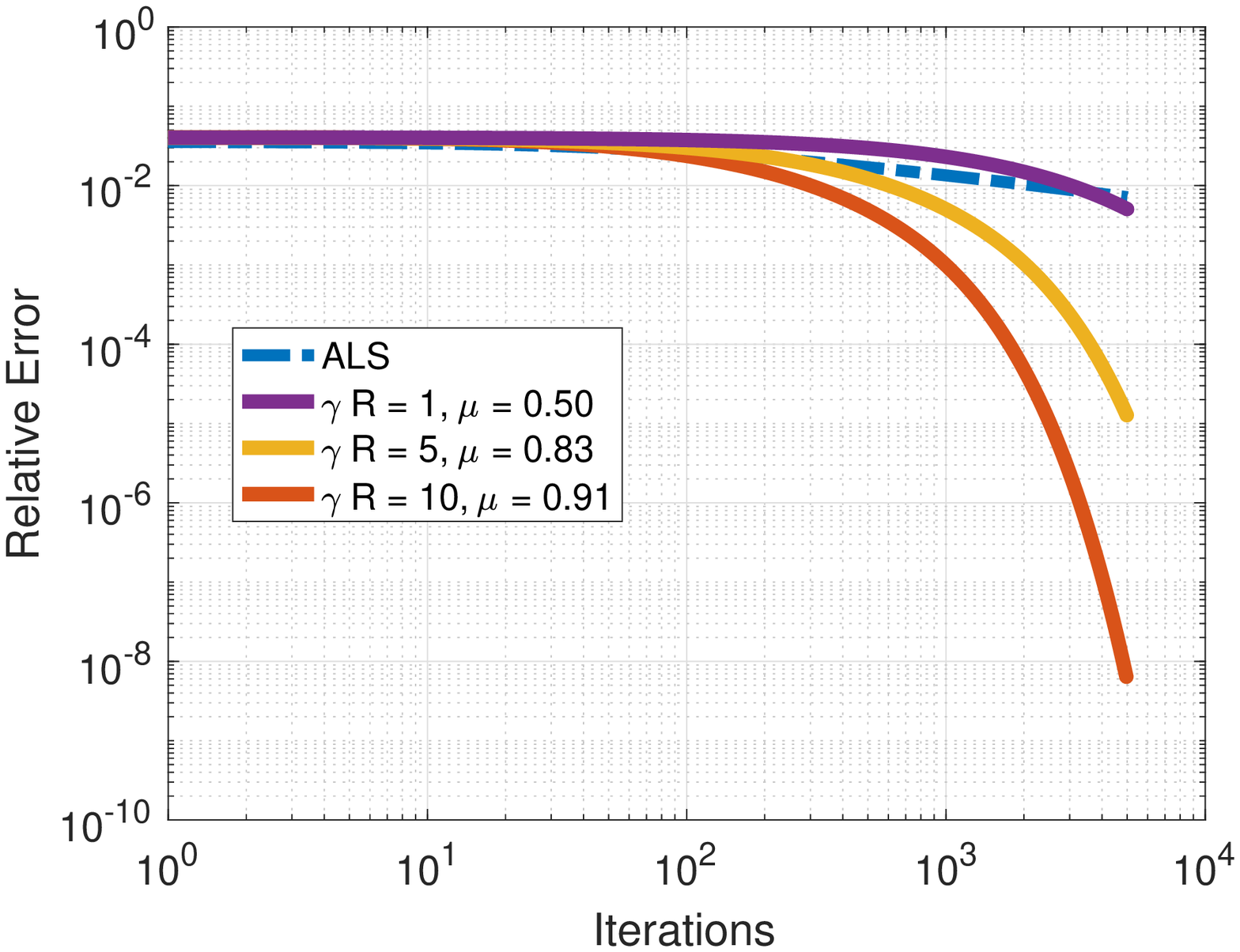}\label{fig_333a2}
}
\subfigure[PARO with regularly adjusted $\gamma$.]{\includegraphics[width=.45\linewidth, trim = 0.0cm 0cm 0cm .0cm,clip=true]{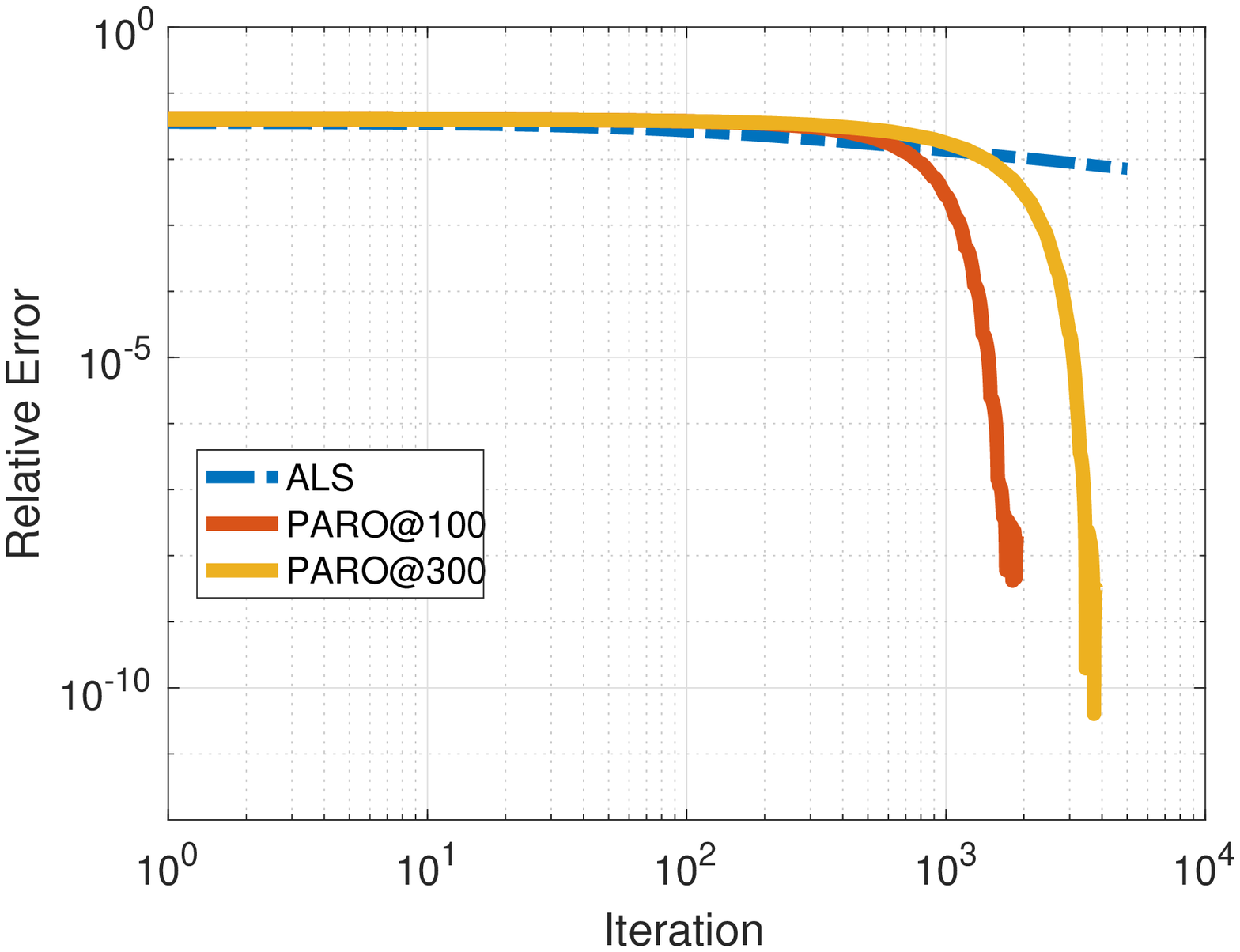}\label{fig_333b2}}
\caption{Performances of ALS and PARO in decomposition of tensors of rank-$R = 23$ associated with the multiplication of two matrices $3\times 3$ and $3 \times 3$.}
\label{fig_3x3x3}
\end{figure}

In this example, we decomposed multiplication tensor associated with multiplication of two matrices of size $3 \times 3$, $3\times 3$. The tensor is of size $9 \times 9 \times 9$ and has rank $R = 23$ \cite{journals/corr/TichavskyPC16}.
We ran PARO with a fixed regularization $\gamma$ which holds $\gamma R = 5$ or 10, or with $\gamma$ adjusted every 100, 200 or 300 iterations.
Success ratios at $10^{-6}$ of PARO with different settings of $\gamma$ are compared with that of ALS in Fig.~\ref{fig_3x3x3}.
The results indicate that PARO with a fixed default $\gamma = \frac{1}{R}$ or $\mu = \frac{1}{2}$ was comparable with the ALS algorithm. The algorithm achieved better performances with higher values of $\gamma$. When $\gamma$ and $\mu$ were regularly adjusted every $K$ iterations, performances of PARO were significantly improved when $K = 200, 300$.
In addition, we illustrate convergence behaviours of PARO and ALS in one run in Fig.~\ref{fig_3x3x3}. ALS and PARO with a fixed $\mu = 1/2$  did not converge in 5000 iterations. However, PARO with $\mu = 0.91$ fully explained the tensor. When $\mu$ was regularly adjusted, PARO converged faster. 
\end{example}

 \begin{example}[Decomposition of tensors with highly collinear loading components]\label{ex_bcd}

\begin{figure}[t]
\centering
\subfigure[]{\includegraphics[width=.48\linewidth, trim = 0.0cm 0cm 0cm .7cm,clip=true]{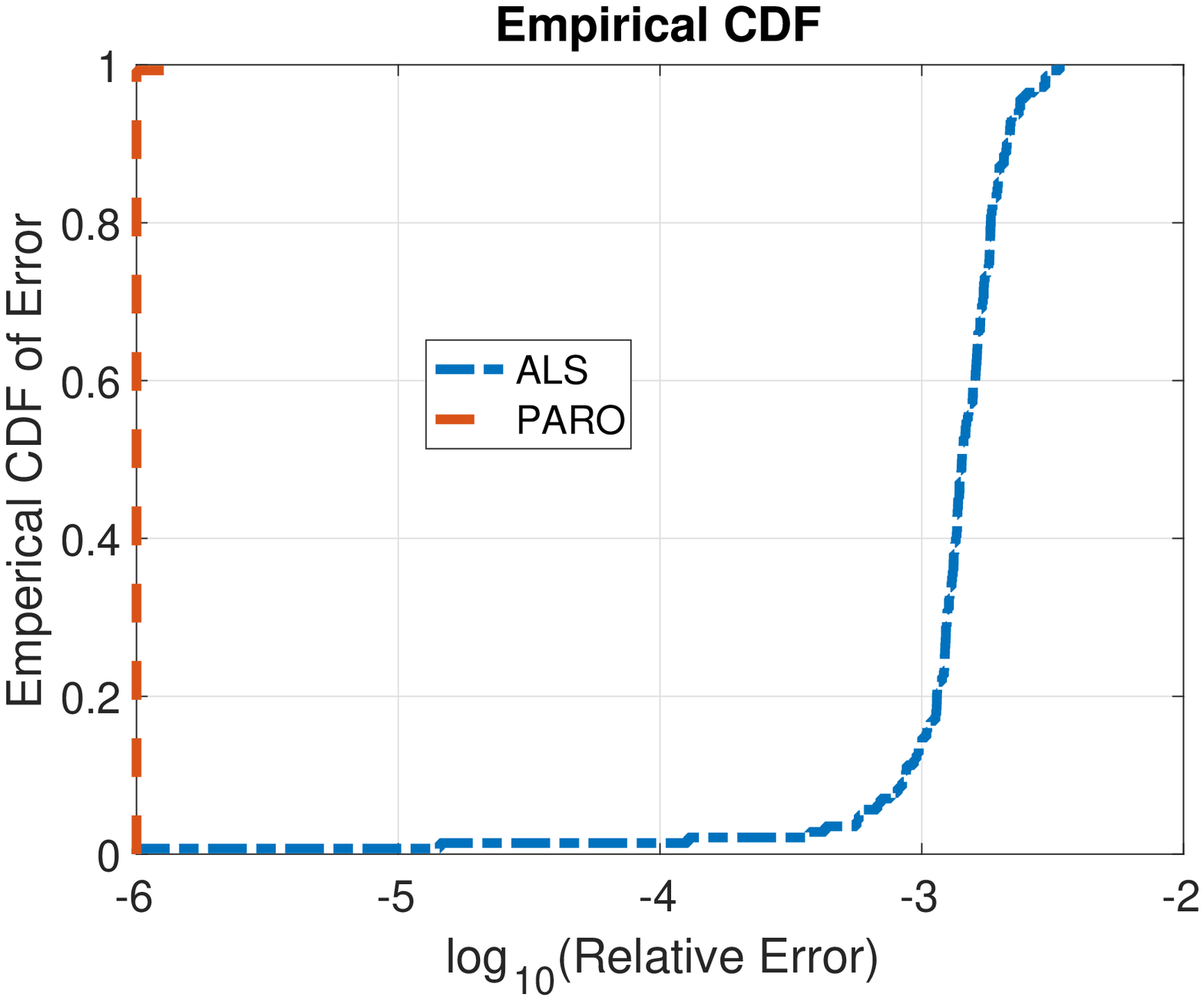}\label{fig_paro_rand0409_b3block_I5N3R8_snrInf_adapt_mu}}
\subfigure[]{\includegraphics[width=.48\linewidth, trim = 0.0cm 0cm 0cm 0cm,clip=true]{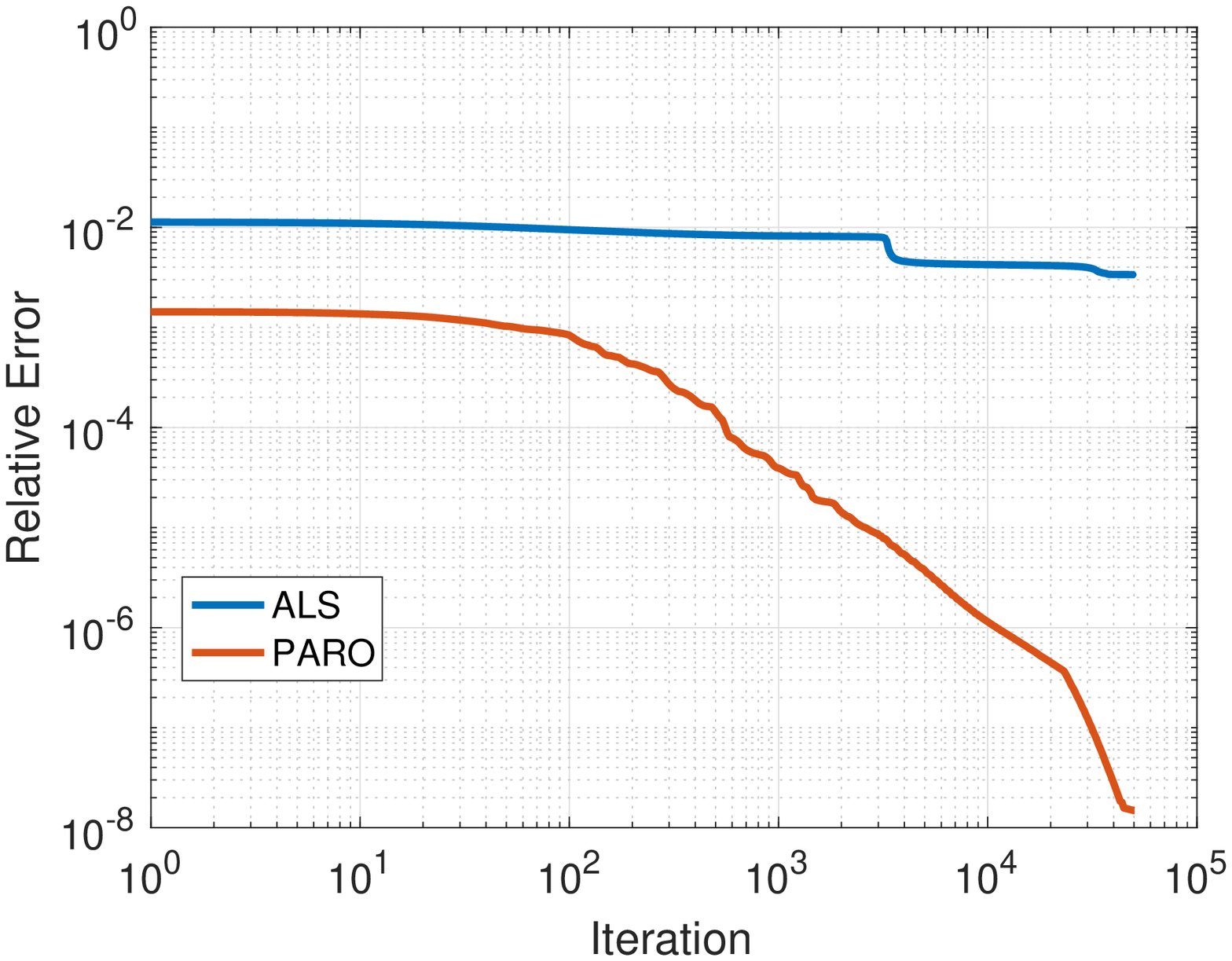}\label{fig_I5R8_paro_vs_als_b3block}}
\caption{(a) Success ratios of the ALS and PARO in Example~\ref{ex_bcd}. (b) Relative errors of the considered algorithms over all runs.}
\label{fig_bcd2blocks}
\end{figure}

 In this example, we decomposed synthesised tensors whose rank exceeded dimensions, and factor matrices comprised highly collinear columns. The tensors were of size $5 \times 5 \times 5$ and rank-8, and composed of two block tensors of rank-4 given in the form of 
 \be
 \tY = \tI \times_1 \bU_{1,1}  \times_2 \bU_{1,2}  \times_3 \bU_{1,3} +  \tI  \times_1 \bU_{2,1}  \times_2 \bU_{2,2}  \times_3 \bU_{2,3} \, , \notag
 \ee
 where $\tI$ represents the diagonal tensor.  
The factor matrices in each block were randomly generated such that the collinearity degree between loading components in a factor matrix was within a range of $[0.95, 0.999]$.
ALS and PARO decomposed  the tensors in 50000 iterations, but could stop earlier if their consecutive approximation errors were different by less than $10^{-9}$.

The ALS almost could not find a CP representation of the noise-free tensors within a relative approximation error range of $10^{-6}$. Its approximation errors were most observed in the interval of $[3.16 \, 10^{-4},  0.0032]$.
The PARO algorithm achieved relative approximation errors smaller than $10^{-6}$ in most of simulation runs. This result is confirmed by the success ratios at $10^{-6}$ of the two algorithms assessed for 100 independent runs, and compared in Fig.~\ref{fig_paro_rand0409_b3block_I5N3R8_snrInf_adapt_mu}.
The relative errors which are most attained by the two algorithms over all runs are illustrated in Fig.~\ref{fig_I5R8_paro_vs_als_b3block}.
 
\end{example}

\section{Conclusions}

In summary, in this paper, we have presented a novel algorithm to update rank-1 tensors in a CPD in parallel. In order to achieve this, we have developed a novel LM algorithm and a novel rotational algorithm for finding best rank-1 tensor approximation. The damped parameter or the step-size in the LM algorithm is optimally determined as root of a degree-$2N$ polynomial, whereas the RORO algorithm seeks a best rank-1 tensor of a $2 \times 2 \times \cdots \times 2$ tensor. Both algorithms have at most the same computational cost as the ALS/HOOI algorithm. The proposed algorithms are verified for decomposition of multiplication tensors and random tensors. 
The PARO algorithm can be implemented in a distributed system with multiple nodes or a system with multiple processing cores in order to decompose tensors of high ranks. 
Finally, the PARO algorithm can be extended to decomposition with nonnegativity constraint, or tensor with incomplete entries. 

\appendices

\section{SVD of a $2\times 2$ Matrix}\label{sec_svd_2x2}

The largest singular value of a $2 \times 2$ matrix $\bW = [w_{i,j}]$ is given by
\be
	\sigma_{\max}^2(\bW) ={\frac{a + \sqrt{a^2 - b^2}}{2}} \label{eq_largest_sv2x2}
\ee
where 
\be
a &=&  \|\bW\|_F^2 = w_{1,1}^2 + w_{1,2}^2 + w_{2,1}^2 + w_{2,2}^2 \, ,\notag \\
b &=&  2 \det(\bW) = 2 (w_{1,1}  w_{2,2} - w_{1,2}  w_{2,1} )  \notag .
\ee

\section{Derivation of Coefficients $a_k$ and $b_k$ in (\ref{eq_opt_rot_eta})}\label{sec::ak_bk}

Derivation of $a(\alpha)$ 
\be
a(\alpha) &=& \|\tW \bar{\times}_3 \boldeta_3 \|_F^2 =  \|\bW_{(1,2)} \boldeta_3\|_2^2 \notag \\&=& \vtr{\bW_{(1,2)}^T \bW_{(1,2)}}^T (\boldeta_3 \otimes \boldeta_3)  \notag \\
&=&  \frac{1}{2}  \left[\cos(2\alpha), \sin(2\alpha), 1\right]   \bQ \vtr{\bW_{(1,2)}^T \bW_{(1,2)}} \,  \notag  \, 
\ee
where $\bW_{(1,2)}$ is mode-(1,2) matricization of $\tW$ and 
\be
\bQ = \left[\begin{array}{rrrr}
    1     & 0 &   0 & -1 \\[-.5em]
     0    & 1 &   1 &  0 \\[-.5em]
     1    & 0 &    0 & 1 
\end{array}\right] \,. \notag 
\ee
This yields the expression of $a_1$, $a_2$ and $a_3$ 
\be
a_1 &=& \frac{1}{2} (\|\bW_1\|_F^2 - \|\bW_2\|_F^2)  \,, \label{eq_a1}\\ 
a_2 &=&   \langle \bW_1, \bW_2 \rangle =  w_ 1 w_ 5 + w_ 2 w_ 6 + w_ 3 w_ 7 + w_ 4  w_ 8 \,, \label{eq_a2}\\
a_3 &=& \frac{1}{2} (\|\bW_1\|_F^2 + \|\bW_2\|_F^2) \,, \label{eq_a3}  
\ee
where $[w_ 1, \ldots, w_ 8] = \vtr{\tW}$. 
Similarly, we can derive expressions of $b_1$, $b_2$ and $b_3$ from 
%
\begin{align}
b(\alpha) &= 2 \det(\tW \bar{\times}_3 \boldeta_3) =  \boldeta_3^T  \bW_{(1,2)}^T \, \bR\,\bW_{(1,2)} \, \boldeta_3  \notag \\
&= \frac{1}{2} \left[\cos(2\alpha), \sin(2\alpha), 1\right]  \, \bQ \, \notag \\
& \quad  (\bw_{1,1,:}  \otimes \bw_{2,2,:}  - \bw_{2,1,:} \otimes \bw_{1,2,:}  - \bw_{1,2,:} \otimes \bw_{2,1,:}  + \bw_{2,2,:} \otimes \bw_{1,1,:}) \notag \, , 
\end{align}
where 
\be
\bR = \left[\begin{array}{rrrr} &&& 1 \\[-.5em] && -1 & \\[-.5em] & -1 & & \\[-.5em] 1 \end{array} \right]\,. \notag 
\ee
Indicating that 
\be
b_1 &=& -w_ 1 w_ 4 + w_ 2  w_ 3  + w_ 5 w_ 8 - w_ 6 w_ 7 \,,\\
b_2 &=&   -w_1  w_8 + w_2 w_7 + w_3 w_6 - w_4 w_5  \,, \label{eq_b2}\\
b_3 &=& -w_ 1 w_ 4 + w_ 2 w_ 3 - w_ 5 w_ 8 + w_ 6 w_ 7 \,.
\ee
We next show that $a_2$ can be zero to simplify the problem after an orthonormal rotation.
Denote by $\bZ$ an orthonormal matrix comprising singular vectors of the mode-3 matricization $\bW_{(3)}$. It is known that rotation of $\tW$ by $\bZ$ along mode-3 yields a tensor ${\tensor{\widetilde{W}}}  = \tW  \, \bar{\times}_3 \, \bZ$ whose two frontal slices are orthogonal, i.e., 
\be
	a_2(\tensor{\widetilde{W}}) = \langle {\tensor{\widetilde{W}}}(:,:,1), {\tensor{\widetilde{W}}}(:,:,2) \rangle = 0\,.
\ee
With this rotation, the two tensors, $\tW$  and ${\tensor{\widetilde{W}}}$, share the largest singular values because 
\be
	\tW  \, \bar{\times}_1 \boldeta_1  \, \bar{\times}_2 \boldeta_2 \, \bar{\times}_3 \boldeta_3 = {\tensor{\widetilde{W}}} \, \bar{\times}_1 \boldeta_1  \, \bar{\times}_2 \boldeta_2 \, \bar{\times}_3  \, \tilde{\boldeta}_3
\ee
where $\tilde{\boldeta}_3 = \bZ \boldeta_3$ remains a unit length vector. 
Hence without loss of generality, we can assume that the frontal slices of the tensor $\tW$  are orthogonal, and $a_2 = 0$.

\section{Derivation of the Polynomial of degree-6 in (\ref{eq_poly_6})}
\label{sec:poly6}

The maximiser $\alpha^{\star}$ to $f(\alpha)$ in (\ref{eq_opt_rot_eta}) is found by setting the derivative $f^{\prime}(\alpha)$ to zero 
\be
 f^{\prime}(\alpha)  
 &=& \frac{a^{\prime}(\alpha)  f(\alpha) - b(\alpha) b^{\prime}(\alpha)}{\sqrt{a^2(\alpha)  - b^2(\alpha)}} = 0 \,,\notag 
\ee
which leads to that 
\be
f(\alpha^{\star}) = \frac{b(\alpha^{\star}) b^{\prime}(\alpha^{\star})}{a^{\prime}(\alpha^{\star})}  \label{eq_opt_rot_eta3}
\ee
and 
\begin{align}
	0 &=  (b(\alpha^{\star}) b^{\prime}(\alpha^{\star}) - a(\alpha^{\star}) \, a^{\prime}(\alpha^{\star}))^2 - {a^{\prime \, 2}(\alpha^{\star})} (a^2(\alpha^{\star})  - b^2(\alpha^{\star}))  \notag \\
	&= b(\alpha^{\star}) (b(\alpha^{\star}) a^{\prime 2}(\alpha^{\star}) + b(\alpha^{\star}) b^{\prime 2}(\alpha^{\star}) - 2 a(\alpha^{\star}) a^{\prime}(\alpha^{\star}) b^{\prime}(\alpha^{\star}) ) \, . \notag
\end{align}
From (\ref{eq_opt_rot_eta3}), since $f(\alpha^{\star}) > 0$, $b(\alpha^{\star})$ must be non-zero, hence, the maximiser $\alpha^{\star}$ is a solution to the following equation 
\be
g(\alpha)  = b(\alpha) a^{\prime 2}(\alpha) + b(\alpha) b^{\prime 2}(\alpha) - 2 a(\alpha) a^{\prime}(\alpha) b^{\prime}(\alpha) = 0 \,.\label{eq_x_222}
\ee
By changing the parameter $\alpha = \arctan(x)$ in the above equation,  and taking into account the equalities 
\be
	\cos(2\arctan(x)) = \frac{1-x^2}{1+x^2}, \quad	\sin(2\arctan(x)) = \frac{2\,x}{1+ x^2},  \notag 
\ee
we come to finding roots of a degree-6 polynomial of coefficients $c_0, c_1, \ldots, c_6$
\be
g(\arctan(x)) = \frac{4}{(1+x^2)^3}  \, (c_6 x^6 + c_5 x^5 +  \cdots  + c_1 x + c_0)  \label{eq_poly_6_2}
\ee
where  
\be
c_6 &=&                                                                           b_2^2 (b_3 - b_1) \label{eq_c6} \\
c_5 &=&                           2  b_2  (2  a_1^2 - 2  a_3  a_1 - 2  b_1^2 + 2  b_3  b_1 + b_2^2) \\
c_4 &=&     4  a_1^2  (b_1 +  b_3) - 8    a_1 a_3 b_1 - 4  b_1^2 (b_1 -    b_3)  + 11  b_1  b_2^2 - b_3  b_2^2 \\
c_3 &=&                                                                         16  b_1^2  b_2 - 4  b_2^3 \\
c_2 &=&  - 4  a_1^2  (b_1  - b_3) - 8  a_1    a_3 b_1 + 4  b_1^2 (b_1 +  b_3)  - 11  b_1  b_2^2 - b_3  b_2^2\\
c_1 &=&                             2  b_2  (2  a_1^2 + 2    a_1  a_3 - 2  b_1^2 - 2    b_1 b_3 + b_2^2)\\
c_0  &=&                                                                         b_2^2  (b_3  + b_1) \,  \label{eq_c0} .
\ee
Note that $a_2 = 0$ and vanishes in the above equations due to the orthonormal rotation in Appendix~\ref{sec::ak_bk}.

 \section{Derivation of Matrices $a(\alpha)$ and $b(\alpha)$ in (\ref{eq_largest_sv2x2_b})}\label{sec:AB_2x2x2x2}
For order-4 tensors, we have 
 \be
a(\alpha) &=& \|\tW \bar{\times}_3 \boldeta_3 \bar{\times}_4 \boldeta_4 \|_F^2 =  \|\bW_{(1,2)} (\boldeta_4  \otimes \boldeta_3)\|_2^2 \notag \\
&=&     \vtr{\bW_{(1,2)}^T \bW_{(1,2)}}^T \bP_{[1 3 2 4]} \,     (\boldeta_4  \otimes \boldeta_4  \otimes \boldeta_3 \otimes  \boldeta_3)   \notag  \,  \\
&=&  \frac{1}{2^2}  \vtr{\bW_{(1,2)}^T \bW_{(1,2)}}^T \bP_{[1 3 2 4]} \, (\bQ  \otimes \bQ) \,  \left( \left[\begin{array}{@{}c@{}} \cos(2\beta) \\[-.5em] \sin(2\beta) \\[-.5em] 1 \end{array}\right]    \otimes \left[\begin{array}{@{}c@{}} \cos(2\alpha) \\[-.5em] \sin(2\alpha) \\[-.5em] 1 \end{array}\right] \right) \notag  \,\\
&=& \left[\begin{array}{@{}c@{}cc} \cos(2\alpha), \sin(2\alpha), 1 \end{array}\right]   \bA \, \left[\begin{array}{@{}c@{}} \cos(2\beta) \\[-.5em] \sin(2\beta) \\[-.5em] 1 \end{array}\right]  \notag 
\ee
where $\bP_{[1 3 2 4]}$ is a commutation matrix which permutes a vectorization of 
a tensor to a vectorization of its mode-$(1, 3, 2, 4)$ permutation, and
\be
\bA &=& \frac{1}{2^2}\bQ^T  (\bW_{1,1,:,:}  \otimes \bW_{1,1,:,:}  + \bW_{2,1,:,:} \otimes \bW_{2,1,:,:}   \notag \\ && \quad + \bW_{1,2,:,:} \otimes \bW_{1,2,:,:}  + \bW_{2,2,:,:} \otimes \bW_{2,2,:,:})  \bQ \,.
\ee
%
%
%
Similarly, we can derive  
\be
b(\alpha) &=& 2 \det(\tW \bar{\times}_3 \boldeta_3 \bar{\times}_4 \boldeta_4) \notag \\ 
&=& 
\frac{1}{2^2}\, \vtr{\bW_{(1,2)}^T  \, \bR \, \bW_{(1,2)}}^T \, \bP_{[1 3 2 4]} \,   (\bQ  \otimes \bQ) \,   
 \left( \left[\begin{array}{@{}c@{}} \cos(2\beta) \\[-.5em] \sin(2\beta) \\[-.5em] 1 \end{array}\right]    \otimes \left[\begin{array}{@{}c@{}} \cos(2\alpha) \\[-.5em] \sin(2\alpha) \\[-.5em] 1 \end{array}\right] \right) \notag \\
&=&  \left[\begin{array}{@{}c@{}cc} \cos(2\alpha), \sin(2\alpha), 1 \end{array}\right]   \bB \, \left[\begin{array}{@{}c@{}} \cos(2\beta) \\[-.5em] \sin(2\beta) \\[-.5em] 1 \end{array}\right]  \notag \, , 
\ee
where  
\be
\bB &=&  \frac{1}{2^2} \bQ^T (\bW_{1,1,:,:}  \otimes \bW_{2,2,:,:}  - \bW_{2,1,:,:} \otimes \bW_{1,2,:,:}   \notag \\ && - \bW_{1,2,:,:} \otimes \bW_{2,1,:,:}  + \bW_{2,2,:,:} \otimes \bW_{1,1,:,:}) \, \bQ \,.
\ee


\section{Derivation of the Two Bi-variate Polynomials of degree-6 in (\ref{eq_bp_41_a})-(\ref{eq_bp_42_a})}\label{sec:bivariate_poly}

Setting gradient of $f(\alpha,\beta)$ w.r.t $\alpha$ and $\beta$ to zeros 
\be
\nabla f(\alpha,\beta) = \frac{f(\alpha,\beta) \nabla a(\alpha,\beta)    -   b(\alpha,\beta)  \, \nabla b(\alpha,\beta)}{\sqrt{a^2(\alpha,\beta) - b^2(\alpha,\beta)}} = 0 \notag 
\ee
gives 
%
\be
{\sqrt{a^2(\alpha,\beta) - b^2(\alpha,\beta)}}  \, \nabla a(\alpha,\beta)    =    b(\alpha,\beta)  \, \nabla b(\alpha,\beta)   -   a(\alpha,\beta) \, \nabla a(\alpha,\beta)     \notag 
\ee  
or the following equations 
\begin{align}
 b(\alpha,\beta) \,   ((a_{\alpha}^{\prime}(\alpha,\beta))^2  + ( b_{\alpha}^{\prime}(\alpha,\beta))^2)   -   2 \, a(\alpha,\beta) \,  a_{\alpha}^{\prime}(\alpha,\beta)    b_{\alpha}^{\prime}(\alpha,\beta)    &= 0  \notag \,,  
 \\
 b(\alpha,\beta) \,   ((a_{\beta}^{\prime}(\alpha,\beta))^2  + ( b_{\beta}^{\prime}(\alpha,\beta))^2)   -   2 \, a(\alpha,\beta) \,  a_{\beta}^{\prime}(\alpha,\beta)    b_{\beta}^{\prime}(\alpha,\beta)    &= 0 \notag \, .
\end{align}
Note that
\be
a_{\alpha}^{\prime}(\alpha,\beta) &=& [\cos(2\alpha), \sin(2\alpha), 1]   \bF^T \, \bA \, [\cos(2\beta), \sin(2\beta), 1]^T \, ,  \notag \\
a_{\beta}^{\prime}(\alpha,\beta) &=& [\cos(2\alpha), \sin(2\alpha), 1]   \, \bA \, \bF \, [\cos(2\beta), \sin(2\beta), 1]^T \, ,  \notag \\
b_{\alpha}^{\prime}(\alpha,\beta) &=& [\cos(2\alpha), \sin(2\alpha), 1]   \bF^T \, \bB \, [\cos(2\beta), \sin(2\beta), 1]^T \, ,  \notag \\
b_{\beta}^{\prime}(\alpha,\beta) &=& [\cos(2\alpha), \sin(2\alpha), 1]   \, \bB \, \bF \, [\cos(2\beta), \sin(2\beta), 1]^T \, ,  \notag 
\ee
where 
\be
\bF = \left[ \begin{array}{rrr}  0 & -1 & 0 \\
1 & 0 & 0 \\
0 & 0 & 0 
\end{array}\right] \, . \notag 
\ee
Next we perform a reparameterization $\alpha = \arctan(x)$ and $\beta = \arctan(z)$. The above two equations become 
\be
g_{\alpha} = \frac{[x^{6}, x^{5}, \ldots, x, 1] \, \bC_1  \,  [z^{6}, z^{5}, \ldots, z, 1]^T }{(1+x^2)^3(1+z^2)^3}   = 0\,,  \label{eq_bp_41}\\
g_{\beta} = \frac{[x^{6}, x^{5}, \ldots, x, 1] \, \bC_2  \,  [z^{6}, z^{5}, \ldots, z, 1]^T}{(1+x^2)^3(1+z^2)^3}   = 0 \,, \label{eq_bp_42} 
\ee
where $\bC_1$ and $\bC_2$ are  two matrices of size $7 \times 7$ defined as 
\begin{align}
\bC_1 &= \bK \, (\bB  \otimes  (\bF^T \bA \otimes \bF^T \bA + \bF^T \bB \otimes \bF^T \bB)  - 2  \bA \otimes \bF^T \bA \otimes \bF^T \, \bB) \bK^T 	\, , \notag \\
\bC_2 &=  \bK \, (\bB  \otimes  (\bA \bF \otimes  \bA  \bF +  \bB \bF\otimes   \bB \bF )  - 2  \bA \otimes \bA \bF\otimes  \bB \bF) \bK^T 
\, ,\notag 
\end{align}
and 
 \be
\bK = \left[ \begin{array}{r@{\hspace{1ex}}r@{\hspace{1ex}}r@{\hspace{1ex}}r@{\hspace{1ex}}r@{\hspace{1ex}}r@{\hspace{1ex}}r@{\hspace{1ex}}r@{\hspace{1ex}}r@{\hspace{1ex}}r@{\hspace{1ex}}r@{\hspace{1ex}}r@{\hspace{1ex}}r@{\hspace{1ex}}r@{\hspace{1ex}}r@{\hspace{1ex}}r@{\hspace{1ex}}r@{\hspace{1ex}}r@{\hspace{1ex}}r@{\hspace{1ex}}r@{\hspace{1ex}}r@{\hspace{1ex}}r@{\hspace{1ex}}r@{\hspace{1ex}}r@{\hspace{1ex}}r@{\hspace{1ex}}r@{\hspace{1ex}}r@{\hspace{1ex}}r@{\hspace{1ex}}r@{\hspace{1ex}}r@{\hspace{1ex}}}
    -1 &    0&     1 &    0&     0 &    0 &    1&     0 &   -1  &   0 &    0   &  0   &  0  &   0  &   0 &    0  &   0   &  0   &  1  &   0  &  -1  &   0   &  0   &  0 &   -1  &   0 &    1\\
     0 &    2&     0 &    2&     0 &   -2 &    0&    -2 &    0  &   2 &    0   & -2   &  0  &   0  &   0 &   -2 &    0  &   2  &   0 &   -2 &    0 &   -2  &   0  &   2&     0 &    2&     0\\
     3 &    0&    -1 &    0&    -4 &    0 &   -1&     0&    -1 &    0&    -4  &   0  &  -4 &    0 &    4&     0 &    4  &   0  &  -1 &    0 &   -1 &    0  &   4  &   0&    -1 &    0&    3\\
     0 &   -4&     0 &   -4&     0 &    0 &    0&     0 &    0  &  -4 &    0   &  0   &  0  &   8  &   0 &    0  &   0   &  4   &  0  &   0  &   0  &   0   &  0   &  4 &    0  &   4 &    0\\
    -3 &    0&    -1 &    0&     4 &    0 &   -1&     0&     1  &   0 &    4   &  0   &  4  &   0  &   4 &    0  &   4   &  0   & -1 &    0 &    1  &   0   &  4   &  0 &    1  &   0 &    3\\
     0 &    2&     0 &    2&     0 &    2 &    0&     2 &    0  &   2 &    0   &  2   &  0  &   0  &   0 &    2  &   0   &  2   &  0  &   2  &   0   &  2    & 0    & 2  &   0   &  2  &   0\\
     1 &    0&     1 &    0&     0 &    0 &    1&     0 &    1  &   0 &    0   &  0   &  0  &   0  &   0 &    0  &   0   &  0   &  1  &   0  &   1   &  0    & 0    & 0  &   1   &  0  &   1\\
     \end{array}
\right] \notag \,.
\ee

\bibliographystyle{IEEEbib}
\bibliography{bibligraphy_thesis,BIBTENSORS2016,BIBTENSORS2017}

\end{document}